\pdfoutput=1
\documentclass[11pt]{article}
\usepackage[T1]{fontenc}
\usepackage[margin=1in]{geometry}
\usepackage{lmodern,microtype,amsmath,amssymb,amsthm,booktabs}
\usepackage{needspace}
\usepackage{capt-of}
\usepackage{enumitem}
\usepackage{tikz}
\usetikzlibrary{decorations.pathreplacing}
\usepackage[hidelinks]{hyperref}
\hypersetup{pdftitle={Bounds for Sidon-Ramsey and van der Waerden numbers},pdfauthor={Edgardo Rold\'an-Pensado, Miguel Raggi and Tania Garibay}}
\newtheorem{theorem}{Theorem}[section]
\newtheorem{lemma}[theorem]{Lemma}
\newtheorem{proposition}[theorem]{Proposition}

\theoremstyle{definition}
\newtheorem{definition}[theorem]{Definition}
\newtheorem{example}[theorem]{Example}

\newcommand{\SR}{\operatorname{SR}}
\newcommand{\SSS}{\mathcal S}
\newcommand{\Z}{\mathbb Z}
\newcommand{\F}{\mathbb F}
\tikzset{
 pt/.style={circle,draw,inner sep=0pt,minimum size=5mm,font=\scriptsize},
 clA/.style={pt,fill=black!80,text=white},
 clB/.style={pt,fill=black!25},
 clC/.style={pt,fill=white},
 cell/.style={draw,minimum size=6mm,inner sep=0pt,font=\scriptsize},
 k1/.style={cell,fill=black!80,text=white},
 k2/.style={cell,fill=black!25},
 k3/.style={cell,fill=white},
}
\newcommand{\WThirteen}{2389}
\newcommand{\WFourteen}{2461}
\newcommand{\WFifteen}{2776}
\newcommand{\WSixteen}{3237}
\newcommand{\WSeventeen}{3602}

\title{Bounds for Sidon--Ramsey and van der Waerden numbers}
\author{Edgardo Rold\'an-Pensado \and Miguel Raggi \and Tania Garibay}
\date{November 2026}

\begin{document}
\maketitle
\begin{abstract}
A set of integers is a Sidon set if all its sums $a+b$ with $a\leq b$ are distinct. The Sidon--Ramsey number $\SR(k)$ is the least $n$ such that $\{1,\ldots,n\}$ cannot be partitioned into $k$ Sidon sets. It is roughly $k^2$, but its exact value is known only for $k\leq7$. We improve the known bounds for $8\leq k\leq20$. For example, $82\leq\SR(8)$, improving $81$, and $103\leq\SR(9)\leq109$, improving $97$ and $111$.

The upper bounds, for $9\leq k\leq12$, rest on a counting observation: a partition of $\{1,\ldots,n\}$ into $k$ Sidon sets must contain several classes of the largest possible size, and there are very few such sets. For eleven and twelve colors not enough of them are pairwise disjoint; for nine and ten colors, integer weights on the points show that the rest of the partition cannot be completed.

The lower bounds are explicit partitions, found by an exact search (eight and nine colors), by a local search that extends known partitions and removes conflicts, or from cyclic constructions. A partition of $\Z_v$ into Sidon sets gives one of $\{1,\ldots,v\}$, and we observe that under a simple condition it can be read slightly past one period. This gives $\SR(q+1)\geq q^2+q+2+m$ for every prime power $q$ and $m(m-1)\leq2q$, and $\SR(p)\geq p^2-p+1+m$ for every prime $p$ and $m(m-1)<2(p-1)$, about $\sqrt{2k}$ more than the classical constructions.

We also improve the lower bounds for the van der Waerden number $W(r,3)$, the least $n$ such that $\{1,\ldots,n\}$ cannot be partitioned into $r$ sets without $3$-term arithmetic progressions, for $13\leq r\leq17$, partly through a product construction. Every result comes with a certificate that a short program checks.
\end{abstract}

\section{Introduction}\label{sec:intro}

\subsection{Sidon sets}\label{sec:intro-sidon}

Write $[n]=\{1,\ldots,n\}$. A set $A$ of integers is a \emph{Sidon set} if all the sums $a+b$ with $a,b\in A$ and $a\leq b$ are different. In other words, in a Sidon set no sum arises in two ways: the value of $a+b$ determines the pair $\{a,b\}$. For example, $\{1,2,5,7\}$ is a Sidon set: its ten sums
\[
 2,\ 3,\ 4,\ 6,\ 7,\ 8,\ 9,\ 10,\ 12,\ 14
\]
($1+1$, $1+2$, $2+2$, $1+5$, $2+5$, $1+7$, $2+7$, $5+5$, $5+7$, $7+7$) are all different. On the other hand, $\{1,2,3\}$ is not a Sidon set, because $1+3=2+2$.

\begin{figure}[ht]
\centering
\begin{tikzpicture}[x=1.2cm,y=0.42cm]
\draw[thick] (0,0) -- (6,0);
\foreach \x in {0,...,6} \draw (\x,0.2) -- (\x,-0.2);
\foreach \x/\l in {0/1,1/2,4/5,6/7} {\fill (\x,0) circle (3pt); \node[above=4pt] at (\x,0) {$\l$};}
\foreach \a/\b/\h/\d in {0/1/1/1,4/6/1/2,1/4/2/3,0/4/3/4,1/6/4/5,0/6/5/6}
 \draw[|-|] (\a,-\h-0.4) -- node[fill=white,inner sep=1.5pt,font=\scriptsize] {$\d$} (\b,-\h-0.4);
\end{tikzpicture}
\caption{The Sidon set $\{1,2,5,7\}$ as a ruler. The six distances between its marks are $1,2,\ldots,6$, each measured exactly once.}
\label{fig:ruler}
\end{figure}
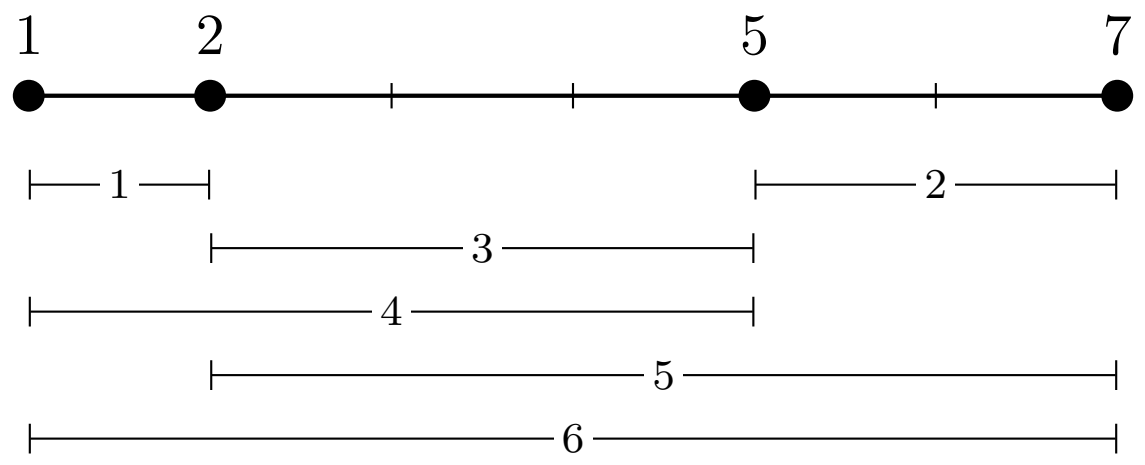

There is an equivalent formulation that is often easier to check: a set is a Sidon set exactly when all the differences $b-a$ with $a<b$ are different, since $a+b=c+d$ is the same as $a-c=d-b$. The differences of $\{1,2,5,7\}$ are $1,2,3,4,5,6$, each once. So a Sidon set is like a ruler with few marks that never measures the same length twice (Figure~\ref{fig:ruler}).

This point of view shows at once that Sidon sets are sparse. A Sidon set with $r$ elements has $\binom r2$ different positive differences, so its largest element exceeds its smallest by at least $\binom r2$. Inside $[n]$ no difference exceeds $n-1$, so a Sidon subset of $[n]$ has fewer than $\sqrt{2n}+1$ elements. A finer argument improves this to about $\sqrt n$.

\subsection{Splitting an interval into Sidon sets}\label{sec:intro-sr}

Since one Sidon set covers only a small part of $[n]$, it is natural to ask how many of them are needed to cover all of it. A \emph{Sidon $k$-partition} of $[n]$ is a partition of $[n]$ into $k$ Sidon sets, which we also call \emph{classes} or \emph{colors}. The \emph{Sidon--Ramsey number} $\SR(k)$ is the least $n$ such that $[n]$ has no Sidon $k$-partition. Equivalently, $\SR(k)$ is the least $n$ such that every coloring of $[n]$ with $k$ colors has a monochromatic solution of $a+b=c+d$ with $\{a,b\}\neq\{c,d\}$; this is the Ramsey-theoretic point of view, and the reason for the name. Removing the largest point of a Sidon $k$-partition of $[n]$ leaves one of $[n-1]$, so $[n]$ has a Sidon $k$-partition exactly when $n<\SR(k)$.

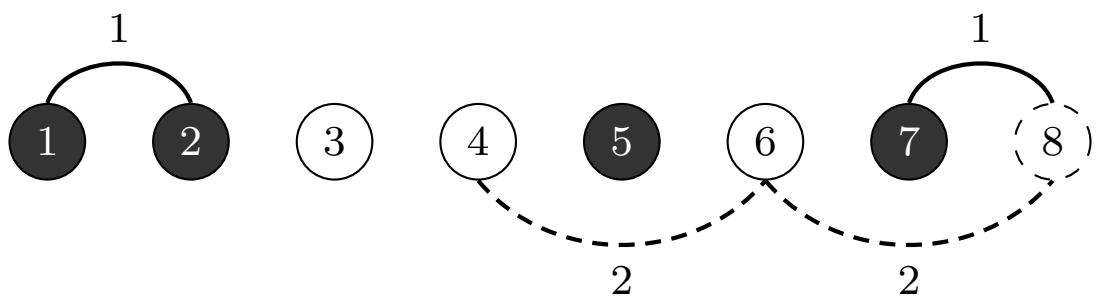
\begin{figure}[ht]
\centering
\begin{tikzpicture}[x=0.95cm]
\foreach \x/\s in {1/clA,2/clA,3/clC,4/clC,5/clA,6/clC,7/clA} \node[\s] (n\x) at (\x,0) {$\x$};
\node[pt,dashed] (n8) at (8,0) {$8$};
\draw[thick] (n1.north) to[bend left=70] node[above,font=\scriptsize] {$1$} (n2.north);
\draw[thick] (n7.north) to[bend left=70] node[above,font=\scriptsize] {$1$} (n8.north);
\draw[thick,densely dashed] (n4.south) to[bend right=50] node[below,font=\scriptsize] {$2$} (n6.south);
\draw[thick,densely dashed] (n6.south) to[bend right=50] node[below,font=\scriptsize] {$2$} (n8.south);
\end{tikzpicture}
\caption{The Sidon $2$-partition $\{1,2,5,7\}\cup\{3,4,6\}$ of $[7]$ (dark and white). The point $8$ fits in neither class: in the dark class it would repeat the difference $1$ ($2-1=8-7$, that is, $1+8=2+7$), and in the white class the difference $2$ ($6-4=8-6$, that is, $4+8=6+6$).}
\label{fig:sr2}
\end{figure}

The smallest cases already show the flavor of the problem. With one color, the question is how long an interval can itself be a Sidon set: $[2]$ is one but $[3]$ is not, so $\SR(1)=3$. With two colors one gets further: $[7]$ splits into $\{1,2,5,7\}\cup\{3,4,6\}$, but $8$ fits in neither class (Figure~\ref{fig:sr2}). That alone does not show that $[8]$ has no Sidon $2$-partition, since the first seven points could be colored differently. Section~\ref{sec:warmup} proves by hand that it has none, so $\SR(2)=8$; that short proof is a model of our upper bounds.

How large is $\SR(k)$? The count above gives a first estimate. If $[n]$ is split into $k$ Sidon sets, each has fewer than $\sqrt{2n}+1$ elements, so $n<k(\sqrt{2n}+1)$, which means roughly $n\leq2k^2$; the finer count improves this to $n\leq k^2+O(k^{3/2})$. In the other direction, the best constructions wrap the interval around a cycle (Section~\ref{sec:cyclic}) and split intervals of length close to $k^2$. So $\SR(k)$ is roughly $k^2$, and near the threshold every class of a Sidon $k$-partition of $[n]$ must be nearly as large as a Sidon subset of $[n]$ can be. The exact values, however, are hard to determine, and they are known only up to seven colors:
\[
 \SR(1),\ldots,\SR(7)=3,\ 8,\ 15,\ 25,\ 36,\ 51,\ 66.
\]
Before this work the best bounds for the next two were $81\leq\SR(8)\leq86$ and $97\leq\SR(9)\leq111$.

\subsection{Results}\label{sec:results}

We improve the lower bounds for every $k$ with $8\leq k\leq20$ and the upper bounds for $9\leq k\leq12$.

\begin{theorem}\label{thm:main}
For $8\leq k\leq20$, $\SR(k)$ lies between the current lower and upper bounds in Table~\ref{tab:bounds}.
\begin{center}\upshape
\begin{tabular}{rr@{\ }lr@{\ }lrr}
\toprule
&\multicolumn{4}{c}{Previous}&\multicolumn{2}{c}{Current}\\
\cmidrule(lr){2-5}\cmidrule(lr){6-7}
$k$&\multicolumn{2}{c}{lower}&\multicolumn{2}{c}{upper}&lower&upper\\
\midrule
8&81&\cite{EP}&86&\cite[Theorem~4.1]{EP}&\textbf{82}&86\\
9&97&\cite{EP}&111&\cite{EP}&\textbf{103}&\textbf{109}\\
10&114&\cite[Table~2]{EP}&141&\cite[Table~2]{EP}&\textbf{119}&\textbf{132}\\
11&133&\cite[Table~2]{EP}&166&\cite[Table~2]{EP}&\textbf{134}&\textbf{158}\\
12&145&\cite[Table~II]{XFL}&193&\cite[Table~2.1]{Dimitromanolakis}&\textbf{155}&\textbf{185}\\
13&170&\cite[Table~II]{XFL}&235&\cite[Table~2.1]{Dimitromanolakis}&\textbf{175}&235\\
14&184&Singer, $q=13$&267&\cite[Table~2.1]{Dimitromanolakis}&\textbf{196}&267\\
15&196&\cite[Table~III]{XFL}&301&\cite[Table~2.1]{Dimitromanolakis}&\textbf{224}&301\\
16&209&\cite[Table~III]{XFL}&321&\cite[Table~2.1]{Dimitromanolakis}&\textbf{251}&321\\
17&274&Singer, $q=16$&392&\cite{OGR24}&\textbf{288}&392\\
18&308&Singer, $q=17$&415&\cite{OGR24}&\textbf{324}&415\\
19&343&Ruzsa, $p=19$&457&\cite{OGR25}&\textbf{352}&457\\
20&382&Singer, $q=19$&521&\cite{OGR27}&\textbf{398}&521\\
\bottomrule
\end{tabular}
\captionof{table}{Bounds for $\SR(k)$. The entries in bold are proved in this work. ``Singer'' and ``Ruzsa'' are the cyclic constructions of Sections~\ref{sec:singer} and~\ref{sec:ruzsa}; Section~\ref{sec:intro-notes} explains how the other previous bounds follow from the cited sources.}
\label{tab:bounds}
\end{center}
\end{theorem}

The upper bounds are proved in Section~\ref{sec:upper}. The proofs are computer-assisted, but checking them needs only integer arithmetic. Each lower bound is $N+1$ for an explicit Sidon $k$-partition of $[N]$ (Section~\ref{sec:search}). All the certificates in this paper, and the programs that check them, are in the verification repository~\cite{SRVerif} at \url{https://github.com/mraggi/sidon-ramsey-verification}. The upper bounds for $k\geq13$ in Table~\ref{tab:bounds} come from Golomb rulers (the ruler bound~\eqref{eq:ruler}) and remain far from the lower bounds.

The classical lower bounds live on a cycle: they split the group $\Z_v$ of residues modulo $v$ into $k$ sets that are Sidon \emph{modulo $v$} (that is, $a+b\equiv c+d\pmod v$ only when $\{a,b\}=\{c,d\}$), and then read the residues as the integers $1,\ldots,v$. Such a partition can be read for longer sequences. Walking around the cycle a second time, when residues meet again, they can create new coincidences, but only between pairs of points exactly one period apart. If the repeated residues happen to lie in different classes, there is no coincidence at all (Lemma~\ref{lem:window}), and a counting argument shows that a suitable stretch of residues always exists. This gives the following general bound.

\begin{theorem}\label{thm:general}
Let $q$ be a prime power and let $m\geq1$ satisfy $m(m-1)\leq2q$. Then
\[
 \SR(q+1)\geq q^2+q+2+m.
\]
Similarly, $\SR(p)\geq p^2-p+1+m$ for every prime $p$ and every $m\geq1$ with $m(m-1)<2(p-1)$.
\end{theorem}

Without the extra stretch, the cyclic construction for $q+1$ colors gives $\SR(q+1)\geq q^2+q+2$; the theorem adds about $\sqrt{2q}$. For small $q$ the longest stretch is much longer than the theorem guarantees, and it gives the bounds for $14$, $17$, $18$, $19$ and $20$ colors in Table~\ref{tab:bounds} (Section~\ref{sec:windows}).

\subsection{How the bounds are proved}\label{sec:ideas}

\paragraph{Upper bounds: the largest classes are rigid.} Near the threshold, a Sidon $k$-partition of $[n]$ is tight: its classes have about $k$ elements each, close to the largest size $M$ that a Sidon subset of $[n]$ can have. A simple count, the \emph{slack} of Section~\ref{sec:slack}, makes this precise: if $kM$ is only slightly larger than $n$, then most classes must have exactly $M$ elements. But Sidon sets of the largest size are rare. For example, $[109]$ has exactly six Sidon subsets with $13$ elements, and no three of them are pairwise disjoint. So a partition would have to contain several disjoint copies of a rare object. For eleven and twelve colors this is already impossible. For nine and ten colors it is not, and we show that the rest of the partition cannot be completed. For this we put integer weights on the points so that every candidate class is heavy while all the points together are light; then not enough candidates can be disjoint (Lemma~\ref{lem:weights}). Good weights are found by linear programming, but checking them needs only additions.

\paragraph{Lower bounds: explicit partitions.} A lower bound $\SR(k)>N$ is proved by exhibiting a Sidon $k$-partition of $[N]$. Ours were found in three ways. For $10\leq k\leq19$, by a local search: start from a partition of a shorter interval, add a point, and recolor points until the coincidences it created (the \emph{conflicts}) disappear. For eight and nine colors, by an exact search: the rigidity used for the upper bounds says that a partition must contain several classes of the largest size, and since there are few of them, we fix them first and complete the partition exhaustively. For twenty colors, by reading a cyclic partition past one period (Theorem~\ref{thm:general}). However a partition was found, checking it is easy: a short program computes the sums in each class.

\subsection{Van der Waerden numbers}\label{sec:intro-vdw}

A Sidon set contains no $3$-term arithmetic progression $a,b,c$, because $a+c=b+b$. Asking only for classes without such progressions relaxes the Sidon condition and leads to van der Waerden numbers. The \emph{van der Waerden number} $W(r,3)$ is the least $n$ such that every coloring of $[n]$ with $r$ colors contains a monochromatic $3$-term arithmetic progression (we write the number of colors first). Hence $\SR(k)\leq W(k,3)$. Sets without progressions can be much denser than Sidon sets, and $W(r,3)$ is much larger than $\SR(r)$: only $W(1,3)=3$, $W(2,3)=9$, $W(3,3)=27$ and $W(4,3)=76$ are known, and for $13\leq r\leq17$ the lower bounds below exceed $10r^2$. We treat the two problems together because they share their methods (cyclic constructions read on an interval, local search that extends colorings and removes conflicts, and certificates checked by short programs), not because of the inequality $\SR(k)\leq W(k,3)$, which is far from tight.

For $5\leq r\leq12$ and for $r=17$ the best published lower bounds come from explicit colorings found by computer, except for $r=6$, which comes from a recursive construction (Table~\ref{tab:vdw-published}). For $13\leq r\leq16$ we found no published coloring. As the previous bounds in Table~\ref{tab:vdw} we take the best bounds that follow from published colorings by the constructions of Section~\ref{sec:vdw-derived}; for $r=17$ this gives $3561$, slightly more than Komkov's $3549$~\cite{Komkov}.

\needspace{18\baselineskip}
\begin{theorem}\label{thm:vdw}
For $13\leq r\leq17$, $W(r,3)>N$, where $N$ is the current lower bound in Table~\ref{tab:vdw}.
\begin{center}\upshape
\begin{minipage}{\linewidth}
\centering
\begin{tabular}{rrr}
\toprule
&\multicolumn{2}{c}{Lower bound}\\
\cmidrule(lr){2-3}
$r$&previous&current\\
\midrule
13&2331&\textbf{\WThirteen}\\
14&2337&\textbf{\WFourteen}\\
15&2676&\textbf{\WFifteen}\\
16&2678&\textbf{\WSixteen}\\
17&3561&\textbf{\WSeventeen}\\
\bottomrule
\end{tabular}
\captionof{table}{Lower bounds $W(r,3)>N$. The previous bounds follow from published colorings (Section~\ref{sec:vdw-derived}); the current ones, in bold, are proved in this work.}
\label{tab:vdw}
\end{minipage}
\end{center}
\end{theorem}

Theorem~\ref{thm:vdw} is proved in Section~\ref{sec:vdw} by explicit colorings, which are too long to print (the longest has $\WSeventeen$ points); the verification repository~\cite{SRVerif} contains them together with a short program that checks them. The colorings for fifteen and sixteen colors start from a product construction, which turns a short coloring with few colors into a much longer one with more colors (Theorem~\ref{thm:palette}); the other colorings, and the last few points of these two, were found by a local search that adds points at the end and removes the monochromatic progressions.

\paragraph{Organization.} Section~\ref{sec:tools} introduces Golomb rulers and two simple facts used throughout. Section~\ref{sec:cyclic} describes the cyclic constructions and how to read them past one period, Section~\ref{sec:upper} proves the upper bounds, Section~\ref{sec:search} describes the searches behind the lower bounds, Section~\ref{sec:vdw} treats van der Waerden numbers, and Section~\ref{sec:open} lists open problems. The appendices contain the integer weights, the complete table of windows, and the computational details: commands, machines, running times, and searches that found nothing. To keep the arguments readable, references are collected in a short ``Notes'' subsection at the end of each section.

\subsection{Notes}\label{sec:intro-notes}

\paragraph{Sidon sets.} Sidon sets go back to Sidon's work on Fourier series~\cite{Sidon} and to Erd\H{o}s and Tur\'an~\cite{ErdosTuran}, who proved that a Sidon subset of $[n]$ has at most $\sqrt n+O(n^{1/4})$ elements; this is the finer count mentioned above. O'Bryant's survey~\cite{OBryant} covers the large literature that followed.

\paragraph{Sidon--Ramsey numbers.} The values $\SR(1),\ldots,\SR(5)$ are due to Liang, Li, Xiu and Xu~\cite{LLXX} (see also Xu, Liang and Luo~\cite{XLL}); $\SR(6)$, $\SR(7)$ and the bounds $81\leq\SR(8)\leq86$ and $97\leq\SR(9)\leq111$ to Espinosa-Garc\'ia and Pellicer~\cite{EP}. For large $k$, Espinosa-Garc\'ia, Montejano, Rold\'an-Pensado and Su\'arez~\cite{EMRS} proved
\[
 k^2-O(k^{c})\leq\SR(k)\leq k^2+Ck^{3/2}+O(k),
\]
where $c\leq1.525$ depends on gaps between primes and $C$ can be taken close to $1.996$. Their lower bound comes from the cyclic construction of Section~\ref{sec:singer}.

\paragraph{The previous bounds in Table~\ref{tab:bounds}.} Through eleven colors they are those collected in~\cite[Table~2]{EP}, together with the two upper bounds proved there. For $k\geq12$ we do not know of a table, so we take the bounds that follow directly from published results. Xiu, Fan and Liang~\cite{XFL} write $H(I,J)=IJ$ when $[IJ]$ can be partitioned into $I$ Sidon sets of size $J$. Their Tables~II and~III give this for $(I,J)=(12,12)$, $(13,13)$, $(15,13)$ and $(16,13)$, hence the lower bounds $145$, $170$, $196$ and $209$. The entries marked Singer and Ruzsa come from the Singer and Ruzsa partitions of Sections~\ref{sec:singer} and~\ref{sec:ruzsa}. The upper bounds come from the ruler bound~\eqref{eq:ruler} with the smallest $M$ such that $g_{M+1}>kM$. The optimal ruler lengths
\[
 g_{17}=199,\quad g_{19}=246,\quad g_{20}=283,\quad g_{21}=333
\]
from~\cite[Table~2.1]{Dimitromanolakis} give $M=16,18,19,20,20$ for $k=12,\ldots,16$, and $g_{24}=425$, established by the distributed computation OGR-24~\cite{OGR24}, gives $M=23$ for $k=17$ and $k=18$. Likewise $g_{25}=480$ and $g_{27}=553$, established by OGR-25~\cite{OGR25} and OGR-27~\cite{OGR27}, give $M=24$ for $k=19$ and $M=26$ for $k=20$. For $k=10,11,12$ the upper bounds obtained in this way, $141$, $166$ and $193$, are improved by Theorem~\ref{thm:main}.

\paragraph{Van der Waerden numbers.} Van der Waerden~\cite{vdW} proved that $W(r,3)$ is finite, and so is its analogue for progressions of any length. The exact values $W(3,3)=27$ and $W(4,3)=76$ are due to Chv\'atal~\cite{Chvatal} and to Beeler and O'Neil~\cite{BeelerONeil}. Except for $r=6$, which follows from a recurrence of Blankenship, Cummings and Taranchuk~\cite{BCT}, the bounds of Table~\ref{tab:vdw-published} come from explicit colorings found by computer: by Heule~\cite{Heule}, with a SAT solver and strong symmetry assumptions, and by Komkov~\cite{Komkov}.

\paragraph{Growth of van der Waerden numbers.} As $r$ grows, $W(r,3)$ grows faster than every power of $r$ but, for large $r$, more slowly than $e^{(\log r)^7}$. Both directions go through $r_3(n)$, the largest size of a progression-free subset of $[n]$. If $r\,r_3(n)<n$, then in every $r$-coloring of $[n]$ some class has more than $r_3(n)$ elements and so contains a progression; hence $W(r,3)\leq n$. Conversely, let $A\subseteq[n]$ be progression-free with $a$ elements. Each $x\in[n]$ lies in exactly $a$ of the $2n-1$ translates $A+t$ with $|t|<n$, so $m$ translates chosen independently at random all miss $x$ with probability $(1-a/(2n-1))^m<e^{-ma/(2n-1)}$. For $m=\lceil(2n-1)\ln n/a\rceil$ this is less than $1/n$, so some $m$ translates cover $[n]$, and coloring each point by a translate that contains it shows $W(m,3)>n$. Behrend's construction~\cite{Behrend} gives $r_3(n)\geq n\,2^{-(2\sqrt2+o(1))\sqrt{\log_2n}}$, and Elsholtz, Hunter, Proske and Sauermann~\cite{EHPS} improved the constant $2\sqrt2$ to $2\sqrt{\log_2(24/7)}\approx2.667$. With the covering argument this gives
\[
 W(r,3)\geq2^{(c-o(1))(\log_2r)^2},\qquad c=\frac1{4\log_2(24/7)}\approx0.141,
\]
where Behrend's constant would give $c=1/8$. This disproves Heule's conjecture~\cite[Conjecture~5]{Heule} that $W(r,\ell)<r^{2\ell}$ for all $r,\ell\geq2$, where $W(r,\ell)$ is defined like $W(r,3)$ with progressions of length $\ell$, for every $\ell\geq3$: a coloring without monochromatic $3$-term progressions has none of length $\ell$, so $W(r,\ell)\geq W(r,3)$, which exceeds $r^{2\ell}$ for all sufficiently large $r$ (even ignoring the $o(1)$, only once $\log_2r>2\ell/c$). In the other direction, Kelley and Meka~\cite{KelleyMeka} proved $r_3(n)\leq n\,e^{-\Omega((\log n)^{1/12})}$, Bloom and Sisask~\cite{BloomSisask} improved the exponent $1/12$ to $1/9$, and Raghavan~\cite{Raghavan} proved $r_3(n)\leq n\,e^{-\Omega((\log n)^{1/6}(\log\log n)^{-1/6})}$. With the first observation, the last bound gives $W(r,3)\leq e^{O((\log r)^6\log\log r)}$. The constants in these upper bounds are not explicit, and the explicit upper bounds we know, such as Gowers' $W(r,3)\leq2^{2^{r^{2^{2^{12}}}}}$~\cite{Gowers}, are far too large to say anything for $r\leq17$; that is why Table~\ref{tab:vdw} gives no upper bounds.

\section{Basic tools}\label{sec:tools}

This section fixes the language of rulers and collects two elementary facts that are used again and again: maps that respect sums carry Sidon sets to Sidon sets, and the classes of a partition have a fixed total amount of unused room.

\subsection{Golomb rulers}\label{sec:rulers}

A Sidon set is a ruler that never measures the same length twice (Section~\ref{sec:intro-sidon} and Figure~\ref{fig:ruler}). A \emph{Golomb ruler} with $r$ marks is a Sidon set $R=\{0=a_1<\cdots<a_r\}$ whose smallest element is $0$; its \emph{length} is $a_r$, and $g_r$ denotes the shortest length of a Golomb ruler with $r$ marks. For instance, $g_2,\ldots,g_6=1,3,6,11,17$. Translation does not change differences, so every Sidon set is a translate of exactly one Golomb ruler (subtract its smallest element). Hence the $r$-element Sidon subsets of $[n]$ are the translates $R+s$, with $1\leq s\leq n-\ell$, of the Golomb rulers $R$ with $r$ marks and length $\ell\leq n-1$. In particular, $[n]$ contains a Sidon set with $r$ elements if and only if $g_r\leq n-1$. Since the $\binom r2$ differences of a ruler are distinct positive integers, $g_r\geq\binom r2$.

\subsection{Maps that respect sums}\label{sec:maps}

We need Sidon sets in groups other than $\Z$. Let $G$ be an abelian group. A set $B\subseteq G$ is a \emph{Sidon set in $G$} if $a+b=c+d$, with $a,b,c,d\in B$, implies $\{a,b\}=\{c,d\}$. For $G=\Z$ this is the definition above; for $G=\Z_v$, the integers modulo $v$, we also say that $B$ is \emph{Sidon modulo $v$}. A subset of a Sidon set is again a Sidon set.

Say that a map $f$ between abelian groups \emph{respects sums} if $a+b=c+d$ implies $f(a)+f(b)=f(c)+f(d)$. Every map $f(x)=\varphi(x)+t$, with $\varphi$ a group homomorphism and $t$ a constant, does, since $f(a)+f(b)=\varphi(a+b)+2t$. If such an $f$ is also injective, it pulls Sidon sets back to Sidon sets: a coincidence $a+b=c+d$ among points that $f$ sends into a Sidon set $B$ becomes the coincidence $f(a)+f(b)=f(c)+f(d)$ in $B$, so $\{f(a),f(b)\}=\{f(c),f(d)\}$, and injectivity gives $\{a,b\}=\{c,d\}$. In other words, a map that respects sums can destroy the Sidon property only by identifying points. We use three instances.
\begin{itemize}
\item \emph{Symmetries of $\Z$.} Translations $x\mapsto x+t$ and reflections $x\mapsto t-x$ are bijections of this form, and so are their inverses, so they map Sidon sets onto Sidon sets. In particular the reflection $x\mapsto n+1-x$ permutes the Sidon subsets of $[n]$.
\item \emph{Symmetries of $\Z_v$.} In the same way, for every unit $u$ modulo $v$ (that is, $\gcd(u,v)=1$) and every $t\in\Z_v$, the affine map $x\mapsto ux+t$ maps sets that are Sidon modulo $v$ onto sets that are Sidon modulo $v$.
\item \emph{From the cycle to the interval.} The map $x\mapsto(x+s-1)\bmod v$ is injective on $[v]$, which consists of $v$ consecutive integers. So if $B$ is Sidon modulo $v$, the integers $x\in[v]$ with $(x+s-1)\bmod v\in B$ form a Sidon set, and a partition of $\Z_v$ into $k$ sets that are Sidon modulo $v$ gives a Sidon $k$-partition of $[v]$; hence $\SR(k)>v$. (With $s=1$: read the residues $1,\ldots,v-1$ as themselves and $0$ as $v$.)
\end{itemize}

The converse of the last item fails: $\{1,2,5\}$ is a Sidon set, but modulo $6$ it has $2+2\equiv5+5$. Being Sidon modulo $v$ is the stronger property, and the cyclic constructions need it.

\subsection{Slack}\label{sec:slack}

\begin{figure}[ht]
\centering
\begin{tikzpicture}[x=0.8cm,y=0.17cm]
\foreach \i/\s in {1/13,2/12,3/12,4/12,5/12,6/12,7/12,8/12,9/12} {
 \fill[black!30] (\i-0.3,0) rectangle (\i+0.3,\s);
 \draw (\i-0.3,0) rectangle (\i+0.3,13);
 \node[below,font=\scriptsize] at (\i,0) {$\s$};
}
\draw[densely dashed] (0.4,13) -- (9.6,13) node[right,font=\scriptsize] {$M=13$};
\end{tikzpicture}
\caption{Slack. Nine classes covering $[109]$, each with room for at most $M=13$ points (the boxes; the shaded part is the size of the class). The shaded parts add up to $109$, so the empty parts add up to $9\cdot13-109=8$ in every such partition. Hence at most eight classes are short of $13$ points, and at least one class has exactly $13$. The figure shows one possible profile of sizes.}
\label{fig:slack}
\end{figure}
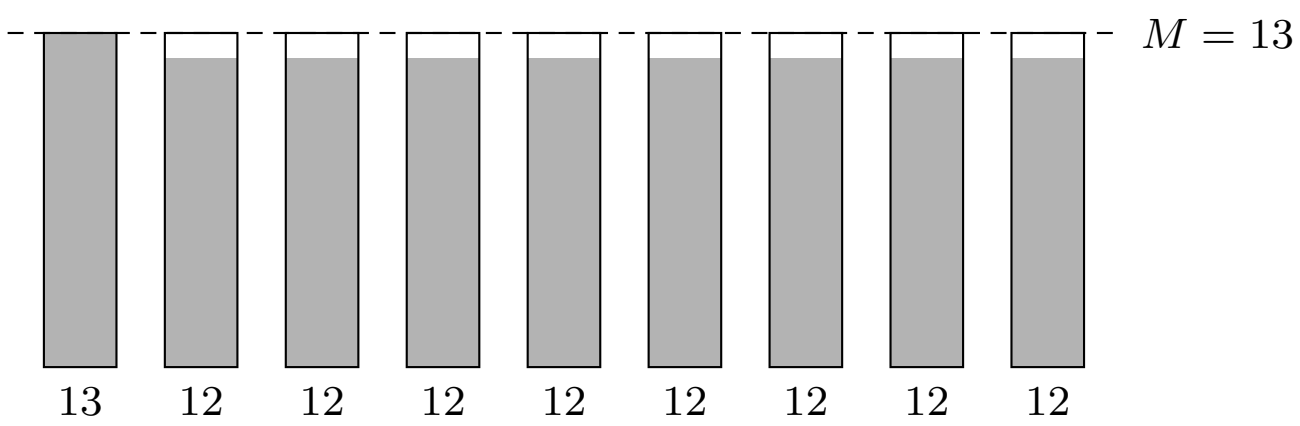

Suppose that $k$ classes must hold $n$ points and that none can hold more than $M$. Think of each class as a box with room for $M$ points. Whatever the partition, the unused room adds up to the same number, the \emph{slack}:
\begin{equation}\label{eq:slack}
 \sum_{C}\bigl(M-|C|\bigr)=kM-n,
\end{equation}
where $C$ runs over the classes. Every term is a nonnegative integer, so $n\leq kM$, and at most $kM-n$ classes have fewer than $M$ points; that is, at least $n-k(M-1)$ classes are full (Figure~\ref{fig:slack}). If the slack is small, almost every box is full.

We use the slack with $M$ the largest size of a Sidon subset of $[n]$, and sometimes a second time for the classes that remain after the largest ones are removed. A first consequence is a quick upper bound, the \emph{ruler bound}:
\begin{equation}\label{eq:ruler}
 \text{if }g_{M+1}>kM,\text{ then }\SR(k)\leq kM+1.
\end{equation}
Indeed, in $[kM+1]$ no difference exceeds $kM<g_{M+1}$, so every Sidon subset has at most $M$ elements, and a Sidon $k$-partition of $[kM+1]$ would have slack $kM-(kM+1)<0$.

\begin{example}\label{ex:slack8}
A Sidon set with $5$ elements has $10$ different positive differences, so it does not fit in $[8]$, where differences are at most $7$. With $k=2$ and $M=4$, the slack of a Sidon $2$-partition of $[8]$ is $2\cdot4-8=0$, so both classes would have exactly four elements. Section~\ref{sec:warmup} takes it from here.
\end{example}

\subsection{Notes}\label{sec:tools-notes}

Golomb rulers have their own literature~\cite{BloomGolomb}. The values of $g_r$ are known for all $r\leq28$; the larger ones were found by distributed computations such as OGR-24~\cite{OGR24}. Golomb rulers have also been used to build disjoint Sidon sets~\cite{Shearer,XFL,XXFL2024}.

\section{Lower bounds from cyclic partitions}\label{sec:cyclic}

A partition of $\Z_v$ into $k$ sets that are Sidon modulo $v$ gives a Sidon $k$-partition of $[v]$ (Section~\ref{sec:maps}). Cyclic groups are easier to work with than intervals because they look the same from every point. In particular, a Sidon set $B$ can be moved around: every translate $B+h$ is again a Sidon set. So if $k$ translates of $B$ cover the group, they give a partition into $k$ Sidon sets at once: put each element in one of the translates that contain it, and each class, being a subset of a Sidon set, is a Sidon set.

Both classical constructions follow this pattern. They differ in the group and in the set $B$.

\subsection{Perfect difference sets}\label{sec:singer}

Consider $D=\{1,2,4\}$ modulo $7$. Its differences are
\[
 2-1=1,\quad 4-2=2,\quad 4-1=3,\quad 1-4\equiv4,\quad 2-4\equiv5,\quad 1-2\equiv6,
\]
so every nonzero residue modulo $7$ is a difference of two elements of $D$ in exactly one way. Read around the cycle, $D$ is a perfect ruler: it measures every distance exactly once.

\begin{definition}
A set $D\subseteq\Z_v$ is a \emph{perfect difference set} if every nonzero residue $x\in\Z_v$ can be written in exactly one way as $x\equiv a-b\pmod v$ with $a,b\in D$.
\end{definition}

If $D$ has $r$ elements, there are $r(r-1)$ differences $a-b$ with $a\neq b$, and they must be exactly the $v-1$ nonzero residues, so $v=r^2-r+1$. With $r=q+1$ this gives $v=q^2+q+1$.

\begin{theorem}[Singer~\cite{Singer}]\label{thm:singer}
For every prime power $q$ there is a perfect difference set with $q+1$ elements modulo $v=q^2+q+1$.
\end{theorem}

The proof uses finite fields: $D$ is a line of the projective plane over the field with $q$ elements (Section~\ref{sec:cyclic-notes}). Here we only need the difference property, and the partition follows from it.

A perfect difference set $D$ is Sidon modulo $v$ for the same reason as in Section~\ref{sec:intro-sidon}: a coincidence $a+b\equiv c+d$ with $\{a,b\}\neq\{c,d\}$ would give two different representations $a-c\equiv d-b$ of one nonzero residue. Its translates $D-d$, $d\in D$, cover $\Z_v$: the residue $0=d-d$ lies in all of them, and a nonzero residue $x\equiv a-d$ lies in $D-d$. Uniqueness of the representation also shows that a nonzero residue lies in only one translate, so the translates share only $0$. Keeping $0$ in one of them gives what we call the \emph{Singer partition} of $\Z_{q^2+q+1}$ into $q+1$ sets that are Sidon modulo $q^2+q+1$: one class with $q+1$ elements and $q$ classes with $q$ each. Hence $\SR(q+1)>q^2+q+1$ for every prime power $q$.

\begin{example}\label{ex:singer2}
For $D=\{1,2,4\}$ modulo $7$ the translates are
\[
 D-1=\{0,1,3\},\qquad D-2=\{0,2,6\},\qquad D-4=\{0,4,5\}.
\]
Keeping $0$ in the first one gives the partition $\{0,1,3\}\cup\{2,6\}\cup\{4,5\}$ of $\Z_7$ and, reading $0$ as $7$, the Sidon $3$-partition $\{1,3,7\}\cup\{2,6\}\cup\{4,5\}$ of $[7]$. For $q=3$, an example of a perfect difference set is $D=\{0,1,3,9\}$ modulo $13$.
\end{example}

\subsection{Exponents and powers}\label{sec:ruzsa}

The second construction gives $p$ classes for a prime $p$, instead of $q+1$. The idea is to pair each exponent $j$ with the power $g^j$. When two such pairs are added, the exponents add, which multiplies the powers, while the powers themselves add. So the sum of two pairs reveals both the sum and the product of the two powers, and two numbers are determined by their sum and their product.

To make this work we need powers that run through all nonzero residues. A \emph{primitive root} modulo $p$ is a number $g$ whose powers
\[
 g^0,\ g^1,\ \ldots,\ g^{p-2}
\]
run through all nonzero residues modulo $p$, each exactly once; one always exists. For example, modulo $5$ the powers of $g=2$ are $1,2,4,3$. Since $g^{p-1}\equiv1\pmod p$, the exponent $j$ only matters modulo $p-1$. We work with pairs $(j,y)$, where $j$ is taken modulo $p-1$ and $y$ modulo $p$, added coordinatewise; that is, in the group $\Z_{p-1}\times\Z_p$. Let
\[
 R=\{(j,g^j):j=0,1,\ldots,p-2\}.
\]
Figure~\ref{fig:ruzsa} shows $R$ for $p=5$ and $g=2$. In the grid of all pairs, $R$ has one cell in each row $j$, and its translate $R+(0,i)$ is $R$ shifted cyclically by $i$ columns. So every cell lies in exactly one of the $p$ translates.

\begin{figure}[ht]
\centering
\begin{tikzpicture}[x=0.75cm,y=0.75cm]
\node[font=\scriptsize,anchor=east] at (-0.6,0.8) {$y=$};
\foreach \y in {0,...,4} \node[font=\scriptsize] at (\y,0.8) {$\y$};
\foreach \j/\g in {0/1,1/2,2/4,3/3} \node[font=\scriptsize,anchor=east] at (-0.6,-\j) {$j=\j$};
\foreach \j/\g in {0/1,1/2,2/4,3/3} \node[font=\scriptsize,anchor=west] at (4.6,-\j) {$g^j=\g$};
\foreach \j/\row in {0/{0,16,12,8,4},1/{5,1,17,13,9},2/{10,6,2,18,14},3/{15,11,7,3,19}}
 \foreach \n [count=\y from 0] in \row \node[k3,minimum size=7mm] at (\y,-\j) {$\n$};
\foreach \j/\y/\n in {0/1/16,1/2/17,2/4/14,3/3/3} \node[k1,minimum size=7mm] at (\y,-\j) {$\n$};
\foreach \j/\y/\n in {0/2/12,1/3/13,2/0/10,3/4/19} \node[k2,minimum size=7mm] at (\y,-\j) {$\n$};
\end{tikzpicture}
\caption{The group $\Z_4\times\Z_5$ for $p=5$ and $g=2$, with row $j$ and column $y$. The dark cells are $R=\{(j,2^j)\}$, and the gray cells are its translate $R+(0,1)$, which is $R$ shifted one column to the right (cyclically). Each cell shows the number $n\in\Z_{20}$ with $n\equiv j\pmod4$ and $n\equiv y\pmod5$; so $R$ corresponds to $\{3,14,16,17\}$ and $R+(0,1)$ to $\{10,12,13,19\}$.}
\label{fig:ruzsa}
\end{figure}

The numbers in Figure~\ref{fig:ruzsa} show how to get back to a cyclic group. Since $p-1$ and $p$ have no common factor, the map $n\mapsto(n\bmod(p-1),\,n\bmod p)$ is a bijection from $\Z_{p(p-1)}$ onto $\Z_{p-1}\times\Z_p$ (this is the Chinese remainder theorem), and it respects sums. So it pulls the $p$ translates of $R$ back to a partition of $\Z_{p(p-1)}$ into $p$ sets, and these are Sidon modulo $p(p-1)$ as soon as $R$ is a Sidon set (Section~\ref{sec:maps}). That is the one point that needs an argument.

\begin{lemma}[Ruzsa~\cite{Ruzsa}]\label{lem:ruzsa}
The set $R$ is a Sidon set in $\Z_{p-1}\times\Z_p$.
\end{lemma}

\begin{proof}
Suppose that $(j_1,x_1)+(j_2,x_2)=(j_3,x_3)+(j_4,x_4)$ with $x_i=g^{j_i}$. The first coordinate says that $j_1+j_2\equiv j_3+j_4\pmod{p-1}$, so the products agree:
\[
 x_1x_2=g^{j_1+j_2}\equiv g^{j_3+j_4}=x_3x_4\pmod p.
\]
The second coordinate says that the sums agree: $x_1+x_2\equiv x_3+x_4\pmod p$. Therefore
\[
 (x_3-x_1)(x_3-x_2)=x_3^2-(x_1+x_2)x_3+x_1x_2\equiv x_3^2-(x_3+x_4)x_3+x_3x_4=0\pmod p.
\]
Since $p$ is prime, $x_3\equiv x_1$ or $x_3\equiv x_2$, and then the sums give $\{x_1,x_2\}=\{x_3,x_4\}$. Each power $x$ comes from a single exponent $j$, so also $\{j_1,j_2\}=\{j_3,j_4\}$.
\end{proof}

We call the resulting partition of $\Z_{p(p-1)}$ into $p$ sets that are Sidon modulo $p(p-1)$ the \emph{Ruzsa partition}. Hence $\SR(p)>p(p-1)$ for every prime $p$.

\begin{example}\label{ex:ruzsa5}
For $p=5$ and $g=2$, Figure~\ref{fig:ruzsa} gives $R=\{3,14,16,17\}$ in $\Z_{20}$. Its sums modulo $20$ are
\[
 6,\ 17,\ 19,\ 0,\ 8,\ 10,\ 11,\ 12,\ 13,\ 14,
\]
all different. Translating by $(0,i)$ means adding the number $n$ with $n\equiv0\pmod4$ and $n\equiv i\pmod5$, that is, adding $0$, $16$, $12$, $8$ or $4$. The translates
\[
 \{3,14,16,17\},\ \{0,1,7,18\},\ \{2,4,5,11\},\ \{6,8,9,15\},\ \{10,12,13,19\}
\]
are $R$, $R+4$, $R+8$, $R+12$ and $R+16$, and they cover $\Z_{20}$ without repetition. Reading $0$ as $20$ gives a Sidon $5$-partition of $[20]$.
\end{example}

\subsection{Reading past one period}\label{sec:windows}

Reading a partition of $\Z_v$ on $[v]$ uses each residue once. Suppose we keep walking around the cycle. On an interval of $v+m$ points, the first $m$ residues are met a second time, so the points $x$ and $x+v$ get the same class for $x\leq m$. This is harmless for a single $x$, but a class that contains two such pairs $\{x,x+v\}$ and $\{y,y+v\}$ has the repeated difference $v$. The lemma below says that this is the only thing that can go wrong.

Let $\Z_v=C_1\cup\cdots\cup C_k$ be a partition into sets that are Sidon modulo $v$. For $s\in\Z$ and $N\geq1$, color $x\in[N]$ by the index $j$ for which $x+s-1$ lies in $C_j$ modulo $v$: we walk around the cycle starting at the residue $s$ and read off the classes. If $N=v$, this is a partition of $[v]$ as in Section~\ref{sec:maps}. If $N=v+m$ with $0<m<v$, the residues $s,s+1,\ldots,s+m-1$ are read twice, by $x$ and by $x+v$, and we call them the \emph{window}.

\begin{figure}[ht]
\centering
\begin{tikzpicture}
\draw[black!25] (0,0) circle (1.25);
\foreach \r/\s in {0/clA,1/clA,2/clB,3/clA,4/clC,5/clC,6/clB} \node[\s] (r\r) at ({90-360/7*\r}:1.25) {$\r$};
\draw[->,thick] ({90-720/7}:1.75) arc[start angle={90-720/7},end angle={90-1440/7},radius=1.75];
\node[font=\scriptsize] at ({90-1080/7}:2.2) {window};
\begin{scope}[xshift=2.6cm,x=0.95cm]
\foreach \x/\r/\s in {1/2/clB,2/3/clA,3/4/clC,4/5/clC,5/6/clB,6/0/clA,7/1/clA,8/2/clB,9/3/clA,10/4/clC} {
 \node[\s] (x\x) at (\x,0) {$\x$};
 \node[font=\scriptsize,below=2pt] at (x\x.south) {$\r$};
}
\draw[decorate,decoration={brace,amplitude=4pt}] (0.65,0.42) -- (3.35,0.42) node[midway,above=5pt,font=\scriptsize] {window};
\draw[decorate,decoration={brace,amplitude=4pt}] (7.65,0.42) -- (10.35,0.42) node[midway,above=5pt,font=\scriptsize] {window again};
\end{scope}
\end{tikzpicture}
\caption{Left: the partition $\{0,1,3\}\cup\{2,6\}\cup\{4,5\}$ of $\Z_7$ from Example~\ref{ex:singer2} (dark, gray, white). Right: reading it from $s=2$ for $10$ steps colors $x\in[10]$ by the class of the residue $x+1$ modulo $7$ (small numbers). The residues $2,3,4$ of the window are read twice; they lie in three different classes, so by Lemma~\ref{lem:window} the result is a Sidon $3$-partition of $[10]$.}
\label{fig:window}
\end{figure}
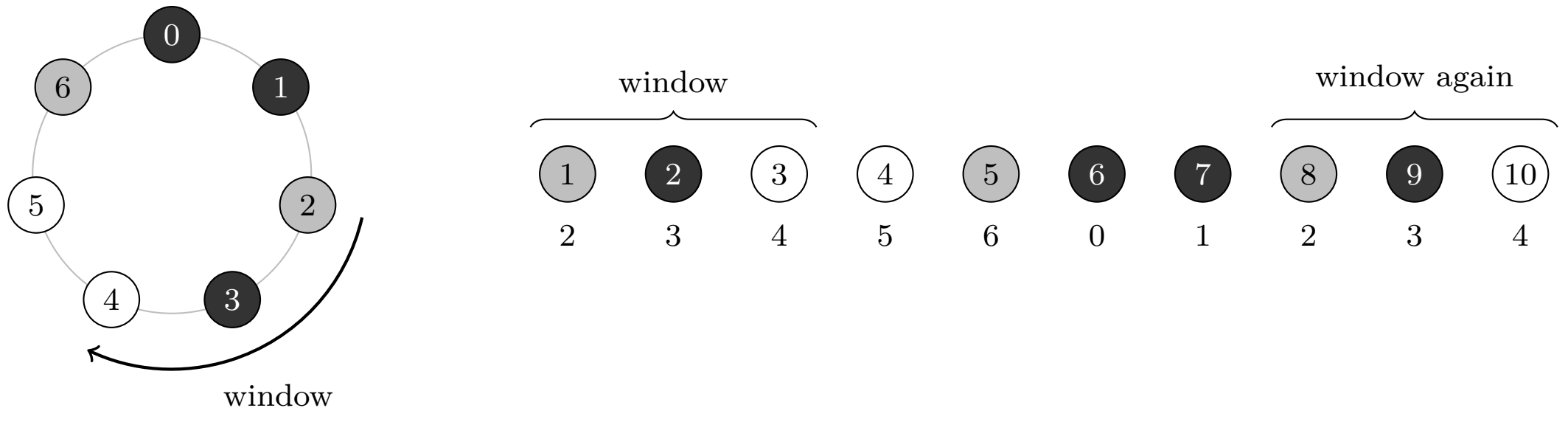

\begin{example}\label{ex:window}
Take the partition $\{0,1,3\}\cup\{2,6\}\cup\{4,5\}$ of $\Z_7$ from Example~\ref{ex:singer2}. The residues $2,3,4$ lie in different classes, so reading from $s=2$ for $10$ steps gives the Sidon $3$-partition
\[
 \{2,6,7,9\}\cup\{1,5,8\}\cup\{3,4,10\}
\]
of $[10]$ (Figure~\ref{fig:window}). Reading from $s=0$ instead fails: the window $0,1,2$ contains $0$ and $1$, which are in the same class, and indeed the class $\{1,2,4,8,9\}$ of $[10]$ contains $2-1=9-8$.
\end{example}

\begin{lemma}\label{lem:window}
If the $m$ residues of the window lie in $m$ different classes, this coloring is a Sidon $k$-partition of $[v+m]$. Hence $\SR(k)>v+m$.
\end{lemma}

\begin{proof}
Let $f(x)$ be the residue of $x+s-1$ modulo $v$, so that the classes are the sets $f^{-1}(C_j)$. The map $f$ respects sums, but it is not injective, so the argument of Section~\ref{sec:maps} does not apply directly; we follow it and see where injectivity is needed. Suppose that $a+b=c+d$ with $a,b,c,d$ in one class $f^{-1}(C_j)$ and $\{a,b\}\neq\{c,d\}$. Then $f(a)+f(b)=f(c)+f(d)$, and since $C_j$ is Sidon modulo $v$, $\{f(a),f(b)\}=\{f(c),f(d)\}$. Exchange $c$ and $d$ if necessary so that $f(a)=f(c)$ and $f(b)=f(d)$. If $a=c$ then also $b=d$, which is excluded; so $a\neq c$. As $a\equiv c\pmod v$ and both lie in $[v+m]$ with $m<v$, we get $|a-c|=v$, and then $b-d=c-a$ shows that $|b-d|=v$ too. So $f(a)$ and $f(b)$ are both residues that are read twice, that is, residues of the window, and both lie in $C_j$. Since different residues of the window lie in different classes, $f(a)=f(b)$. But then $a$, $b$, $c$, $d$ all lie in $\{x,x+v\}$ for the same $x$, with $a\neq c$ and $b\neq d$; from $a-c=d-b$ we get $\{a,b\}=\{x,x+v\}=\{c,d\}$, a contradiction.
\end{proof}

\paragraph{Many partitions from one.} A single partition of $\Z_v$ gives many others. For every unit $u$ modulo $v$ and every $t\in\Z_v$, the map $x\mapsto ux+t$ preserves the Sidon property (Section~\ref{sec:maps}), so the sets $u^{-1}(C_j-t)=\{y:uy+t\in C_j\}$ again form a partition into sets that are Sidon modulo $v$. Reading this partition from $s=1$ colors $x$ by the class of $ux+t$ in the original partition, and the window consists of the residues $u+t,2u+t,\ldots,mu+t$. There are $v\varphi(v)$ such maps, where $\varphi(v)$ is the number of units modulo $v$ ($39312$, $93942$ and $96012$ for the Singer partitions with $q=16,17,19$); they need not give pairwise different partitions of $[v]$, since a multiplier of the difference set that fixes the class containing $0$ maps the partition to itself. For Singer partitions the residue $0$ lies on all the translates $D-d$, so we put it in a class not used by the other residues of the window. Finding the longest window over all images is a quick finite computation; Table~\ref{tab:windows} shows the results for the constructions we use with at most $20$ colors.

\begin{table}[ht]
\centering
\begin{tabular}{lrrrrl}
\toprule
Construction&$v$&$k$&$m$&$(u,t)$&Lower bound\\
\midrule
Singer, $q=13$&183&14&12&$(17,47)$&$\SR(14)\geq196$\\
Ruzsa, $p=17$&272&17&8&$(1,9)$&$\SR(17)\geq281$\\
Singer, $q=16$&273&17&14&$(19,188)$&$\SR(17)\geq288$\\
Singer, $q=17$&307&18&16&$(9,293)$&$\SR(18)\geq324$\\
Ruzsa, $p=19$&342&19&9&$(1,9)$&$\SR(19)\geq352$\\
Singer, $q=19$&381&20&16&$(44,147)$&$\SR(20)\geq398$\\
\bottomrule
\end{tabular}
\caption{The longest windows over all affine images $x\mapsto ux+t$; the pair $(u,t)$ is one image that attains it, and $m$ is the window length. ``Singer'' and ``Ruzsa'' are the partitions of Sections~\ref{sec:singer} and~\ref{sec:ruzsa}.}
\label{tab:windows}
\end{table}

For $14$, $17$, $18$ and $19$ colors these bounds equal those reached by the local search of Section~\ref{sec:search}; for $20$ colors the window gives one more point ($[397]$ instead of $[396]$). The partitions of $[287]$, $[323]$ and $[396]$ found by the search, started from affine images of Singer partitions, are almost periodic: at $44$ of the $45$ points $v+x$ past the first period the color equals that of $x$.

\paragraph{A window always exists.} Why should a good window exist at all? A window is bad when it contains two residues $x$ and $x+\delta$ of the same class. Each class contains at most one such pair for each distance $\delta$, so bad pairs are rare, and averaging over the $v$ possible starting points shows that some window of moderate length avoids them all.

\begin{proposition}\label{prop:average}
Let $\Z_v=C_1\cup\cdots\cup C_k$ be a partition into sets that are Sidon modulo $v$, and let $1\leq m<v$ with $k\,m(m-1)<2v$. Then some window of length $m$ has its residues in different classes, so $\SR(k)>v+m$.
\end{proposition}

\begin{proof}
Fix $\delta$ with $1\leq\delta<m$. A class cannot contain two pairs $\{a,a+\delta\}$ and $\{b,b+\delta\}$ with $a\neq b$, because $a+(b+\delta)\equiv b+(a+\delta)$ would be a nontrivial solution. So at most $k$ residues $x$ have $x$ and $x+\delta$ in the same class. Each such pair lies in $m-\delta$ of the $v$ windows $\{s,\ldots,s+m-1\}$. Summing over $\delta$, the $v$ windows contain at most $k\sum_{\delta=1}^{m-1}(m-\delta)=k\,m(m-1)/2<v$ such pairs in total, so some window contains none, and Lemma~\ref{lem:window} applies.
\end{proof}

Theorem~\ref{thm:general} is Proposition~\ref{prop:average} applied to the Singer partition, with $k=q+1$ and $v=q^2+q+1$ (where $m(m-1)\leq2q$ gives $k\,m(m-1)\leq2q^2+2q<2v$), and to the Ruzsa partition, with $k=p$ and $v=p(p-1)$ (where the condition is $m(m-1)<2(p-1)$).

How long are the best windows? Table~\ref{tab:windows-large} in Appendix~\ref{app:windows} lists them for all Singer partitions with $q\leq64$ and all Ruzsa partitions with $p\leq67$. For $q\leq8$ the window has the largest possible length $k=q+1$, so $\SR(q+1)\geq(q+1)^2+2$ for these $q$ (all far below the true values). After that the ratio $m/k$ slowly decreases, from about $0.9$ near $k=18$ to $0.58$ for $q=64$. A heuristic explains the trend: $m$ residues fall into different classes with probability about $\prod_{i<m}(1-i/k)\approx e^{-m^2/2k}$, and there are about $v^2\approx k^4$ windows, which suggests that the longest window has length of order $\sqrt{k\log k}$. Since the residues of a window lie in different classes, $m\leq k$, so windows affect neither the leading term $k^2$ nor the error terms of the asymptotic bounds in Section~\ref{sec:intro-notes}.

\subsection{Notes}\label{sec:cyclic-notes}

Perfect difference sets are due to Singer~\cite{Singer}, who obtained them from the projective plane over the finite field $\F_q$; Ruzsa's construction is from~\cite{Ruzsa}. Partitions of cyclic groups into Sidon sets are also studied in~\cite{DMRT}. We believe that Lemma~\ref{lem:window} is probably known, but we could not find it in the literature.

\section{Upper bounds}\label{sec:upper}

All our upper bounds follow one plan.
\begin{enumerate}
\item \emph{Slack.} By the slack count~\eqref{eq:slack}, a Sidon $k$-partition of $[n]$ must contain several classes of the largest possible size $M$.
\item \emph{Rigidity.} Sidon subsets of $[n]$ with $M$ elements are rare, so these classes can be listed.
\item \emph{Completion.} The listed sets cannot be combined into a partition: either not enough of them are pairwise disjoint (eleven and twelve colors), or what is left cannot be split as it would have to be (nine and ten colors), which we prove with weights.
\end{enumerate}
We first carry out the plan by hand for two colors (Section~\ref{sec:warmup}), then for eleven and twelve colors, where step~3 is easiest. Section~\ref{sec:weights} introduces the weights, and Sections~\ref{sec:nine} and~\ref{sec:ten} prove the bounds for nine and ten colors. Section~\ref{sec:stop} explains why weights alone cannot go below $109$ for nine colors, and Section~\ref{sec:computation} describes the computations.

\subsection{A warm-up: two colors}\label{sec:warmup}

\begin{proposition}\label{prop:sr2}
$[8]$ has no Sidon $2$-partition. Hence $\SR(2)=8$.
\end{proposition}

\begin{proof}
\emph{Slack.} By Example~\ref{ex:slack8}, both classes of a Sidon $2$-partition of $[8]$ would have exactly four elements.

\emph{Rigidity.} One of the two classes, call it $B$, does not contain $1$, so $B\subseteq\{2,\ldots,8\}$. Its six differences are different numbers between $1$ and $6$, so they are exactly $1,\ldots,6$. In particular $B$ contains $2$ and $8$, and $B-2=\{0,a,b,6\}$ is a Golomb ruler that measures each length from $1$ to $6$ once. Checking the ten choices $1\leq a<b\leq5$ leaves only $\{0,1,4,6\}$ and its mirror image $\{0,2,5,6\}$, so $B=\{2,3,6,8\}$ or $B=\{2,4,7,8\}$.

\emph{Completion.} The other class would be the complement of $B$, that is, $\{1,4,5,7\}$ or $\{1,3,5,6\}$. Neither is a Sidon set: $1+7=4+4$ and $1+5=3+3$.

So $\SR(2)\leq8$, and the partition of $[7]$ in Figure~\ref{fig:sr2} gives $\SR(2)>7$.
\end{proof}

\subsection{Eleven and twelve colors}\label{sec:upper1112}

When $n-1$ is only slightly larger than the ruler length $g_M$, the $M$-element Sidon subsets of $[n]$ are few: they are the translates $R+s$, with $1\leq s\leq n-\ell$, of the Golomb rulers $R$ with $M$ marks and length $\ell\leq n-1$ (Section~\ref{sec:rulers}). Table~\ref{tab:near} counts the rulers we need. For eleven and twelve colors, step~3 of the plan only asks how many of these translates can be pairwise disjoint.

\begin{table}[ht]
\centering
\begin{tabular}{lrrrrrrrr}
\toprule
$M=14$, length&$127$&$128$&$129$&$130$&$131$&&&\\
rulers&$2$&$2$&$6$&$2$&$12$&&&\\
\midrule
$M=15$, length&$151$&$152$&$153$&$154$&$155$&$156$&$157$&\\
rulers&$2$&$4$&$4$&$0$&$10$&$14$&$24$&\\
\midrule
$M=16$, length&$177$&$178$&$179$&$180$&$181$&$182$&$183$&$184$\\
rulers&$2$&$2$&$6$&$6$&$12$&$14$&$36$&$42$\\
\bottomrule
\end{tabular}
\caption{The number of Golomb rulers with $M$ marks and the given length; a ruler and its reflection are counted separately.}
\label{tab:near}
\end{table}

\begin{proof}[Proof that $\SR(11)\leq158$ and $\SR(12)\leq185$]
\emph{Eleven colors.} Since $g_{16}=177>157$, every Sidon subset of $[158]$ has at most $15$ elements. The slack of a Sidon $11$-partition of $[158]$ would be $11\cdot15-158=7$, so at least $11-7=4$ of its classes, which are pairwise disjoint, have $15$ elements. The $15$-element Sidon subsets of $[158]$ are the $140$ translates of the $58$ rulers of Table~\ref{tab:near} with $15$ marks, and no four of them are pairwise disjoint ($94$ triples are).

\emph{Twelve colors.} Since $g_{17}=199>184$, every Sidon subset of $[185]$ has at most $16$ elements. The slack of a Sidon $12$-partition of $[185]$ would be $12\cdot16-185=7$, so at least $5$ pairwise disjoint classes have $16$ elements. The $16$-element Sidon subsets of $[185]$ are the $300$ translates of the $120$ rulers of Table~\ref{tab:near} with $16$ marks, and no five of them are pairwise disjoint ($37$ quadruples are).
\end{proof}

\subsection{Weights}\label{sec:weights}

For nine and ten colors there are enough disjoint classes of the largest size, and we must show that the rest of the partition cannot be completed. The tool is a price list. Give every point $x$ of the set $U$ to be split a price $w(x)$, and let a set cost the sum of the prices of its points. If $t$ disjoint candidates $S_1,\ldots,S_t$ leave the points $T$ uncovered, the price of all the points splits as
\[
 w(U)=w(S_1)+\cdots+w(S_t)+w(T).
\]
If the prices are nonnegative and every candidate costs at least $L$, the right side is at least $tL$. So if all the points together cost less than $tL$, no $t$ candidates are disjoint.

The same bookkeeping is sharper when the candidates must cover all of $U$ except exactly $f$ points, as happens when the remaining classes have prescribed sizes. Then $T$ has exactly $f$ points, so it costs at least the sum of the $f$ smallest prices, not merely $0$. This bound holds whatever the signs of the prices, so prices may now be negative, and it is enough that every candidate costs at least $0$. If all the points together cost less than the sum of the $f$ smallest prices, the candidates cannot cover all but $f$ points.

\begin{lemma}[Weights]\label{lem:weights}
Let $\mathcal F$ be a family of subsets of a finite set $U$, let $w:U\to\Z$, and write $w(S)=\sum_{x\in S}w(x)$.
\begin{enumerate}[label=(\alph*)]
\item\label{it:pack} If $w\geq0$ and $w(S)\geq L>0$ for every $S\in\mathcal F$, then at most $w(U)/L$ members of $\mathcal F$ are pairwise disjoint.
\item\label{it:cover} If $w(S)\geq0$ for every $S\in\mathcal F$, and $w(U)$ is smaller than the sum of the $f$ smallest values of $w$, then $U$ is not the union of pairwise disjoint members of $\mathcal F$ and a set of $f$ further points.
\end{enumerate}
\end{lemma}

For $f=0$, condition~\ref{it:cover} simply says $w(U)<0$.

\begin{example}\label{ex:weights8}
Here is a second proof of the completion step of Proposition~\ref{prop:sr2}. Give the points $1,\ldots,8$ the weights
\[
 3,\ 2,\ 1,\ 0,\ 0,\ 1,\ 2,\ 3.
\]
The ten $4$-element Sidon subsets of $[8]$ and their weights are
\[
\begin{array}{lllll}
\{1,2,4,8\}:8,&\{1,2,5,7\}:7,&\{1,2,6,8\}:9,&\{1,3,4,8\}:7,&\{1,3,6,7\}:7,\\
\{1,3,7,8\}:9,&\{1,5,6,8\}:7,&\{1,5,7,8\}:8,&\{2,3,6,8\}:7,&\{2,4,7,8\}:7.
\end{array}
\]
Every one has weight at least $7$, while all eight points weigh $12<2\cdot7$. By Lemma~\ref{lem:weights}\ref{it:pack}, no two of these sets are disjoint, so $[8]$ has no Sidon $2$-partition. The weights reflect the fact that Sidon sets must be spread out: a $4$-element Sidon subset of $[8]$ has its largest element at least $6$ above its smallest, so it uses points near the ends, and those are the expensive ones.
\end{example}

\begin{example}\label{ex:signed8}
Part~\ref{it:cover} reaches the same conclusion with weights of both signs. Multiply the weights of Example~\ref{ex:weights8} by $4$ and subtract $7$:
\[
 5,\ 1,\ -3,\ -7,\ -7,\ -3,\ 1,\ 5.
\]
A $4$-element Sidon subset of $[8]$ now weighs at least $4\cdot7-4\cdot7=0$, while all eight points weigh $4\cdot12-8\cdot7=-8<0$. By Lemma~\ref{lem:weights}\ref{it:cover} with $f=0$, $[8]$ is not a union of disjoint $4$-element Sidon sets. When $f=0$ and all candidates have the same size, such a shift turns weights for either part into weights for the other, so the two parts prove the same things. Part~\ref{it:cover} is stronger when $f>0$, because it charges the leftover points; the proof for ten colors (Section~\ref{sec:ten}) uses it with a class of $12$ points left over.
\end{example}

\paragraph{Where weights come from.} The best weights for Lemma~\ref{lem:weights}\ref{it:pack} solve a linear program: minimize the total weight $\sum_xy_x$ over real numbers $y_x\geq0$ that give every member of $\mathcal F$ weight at least $1$. If its optimum $\tau^*$ is less than $t$, multiplying an optimal solution by a large integer $L$ and rounding up gives integer weights that exclude $t$ pairwise disjoint members; only these integers enter a proof, so the accuracy of the solver does not matter. If $\tau^*\geq t$, no weights can do it. Linear programming duality identifies $\tau^*$: it is the largest total $\sum_Sy_S$ of numbers $y_S\geq0$, one for each member of $\mathcal F$, such that every point lies in members of total value at most $1$. Such $y$ is a \emph{fractional packing}. For instance, the weights of Example~\ref{ex:weights8} divided by $7$ give every member weight at least $1$ and total $12/7$, so $\tau^*\leq12/7<2$ there. In fact $\tau^*=12/7$: the values $\frac27$, $\frac17$, $\frac37$, $\frac17$, $\frac37$, $\frac27$ on the sets $\{1,2,5,7\}$, $\{1,3,4,8\}$, $\{1,3,6,7\}$, $\{1,5,6,8\}$, $\{2,3,6,8\}$, $\{2,4,7,8\}$ form a fractional packing with this total.

We never need duality in a proof, only its easy direction, which explains when weights must fail. If numbers $y_S\geq0$ cover every point exactly once (a \emph{fractional tiling}, with every point in members of total value exactly $1$), then for any $w$ that gives every member weight at least $L$,
\begin{equation}\label{eq:tiling}
 w(U)=\sum_{S\in\mathcal F}y_S\,w(S)\geq L\sum_Sy_S,
\end{equation}
so Lemma~\ref{lem:weights}\ref{it:pack} cannot exclude $\sum_Sy_S$ disjoint members, even if there are none.

\begin{figure}[ht]
\centering
\begin{tikzpicture}[x=0.55cm,y=0.45cm]
\foreach \x in {1,...,15} \node[font=\scriptsize] at (\x,0.2) {$\x$};
\foreach \row/\set in {1/{1,5,10,11,13},2/{3,5,6,11,15},3/{1,8,9,12,14},4/{2,4,7,8,15},5/{2,4,9,10,13},6/{3,6,7,12,14}} {
 \node[font=\scriptsize,anchor=east] at (0.3,-\row) {$S_\row$};
 \foreach \x in {1,...,15} \fill[black!15] (\x,-\row) circle (2pt);
 \foreach \x in \set \fill (\x,-\row) circle (4pt);
}
\draw[black!40] (0.5,-6.6) -- (15.5,-6.6);
\node[font=\scriptsize,anchor=east] at (0.3,-7.2) {count};
\foreach \x in {1,...,15} \node[font=\scriptsize] at (\x,-7.2) {$2$};
\end{tikzpicture}
\caption{Six $5$-element Sidon subsets of $[15]$, one per row; $S_1$ and $S_2$, $S_3$ and $S_4$, $S_5$ and $S_6$ are mirror images under $x\mapsto16-x$. Every point lies in exactly two of them. With the value $\frac12$ on each set, this is a fractional tiling of $[15]$ with total $3$, although no three $5$-element Sidon subsets of $[15]$ are pairwise disjoint.}
\label{fig:tiling15}
\end{figure}
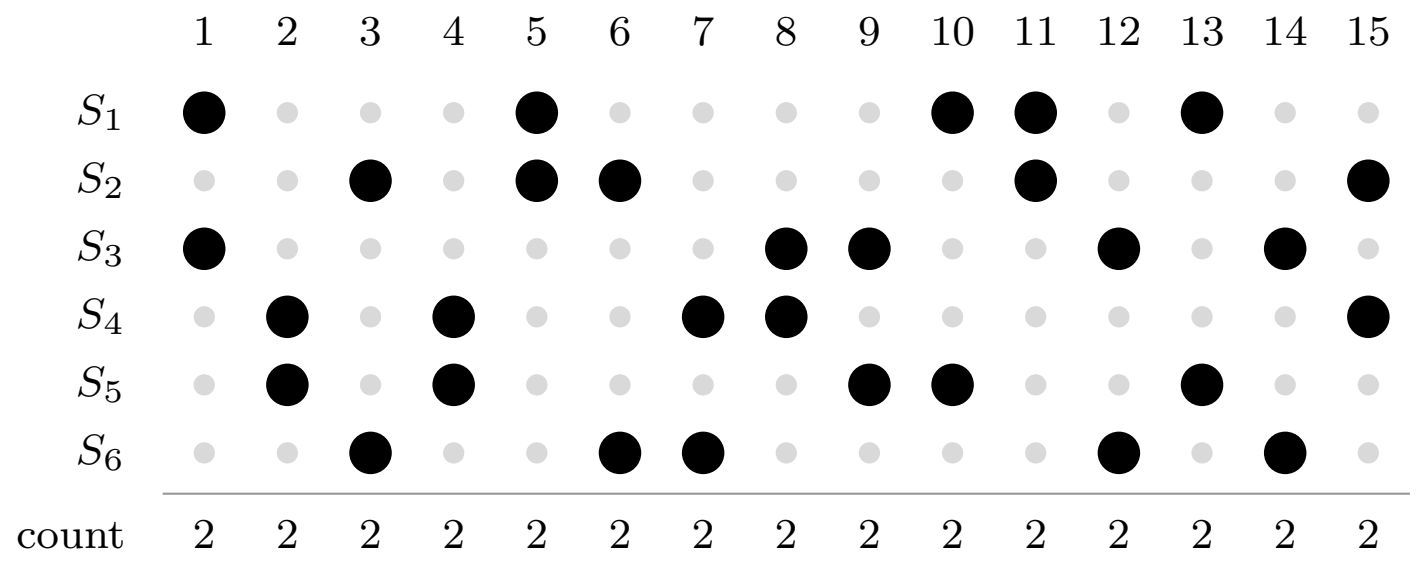

\begin{example}\label{ex:tiling15}
A Sidon set with $6$ elements has $15$ different positive differences and does not fit in $[15]$, so a Sidon $3$-partition of $[15]$ would have slack $3\cdot5-15=0$, that is, three pairwise disjoint $5$-element Sidon sets. Since $\SR(3)=15$, there are no such three among the $156$ five-element Sidon subsets of $[15]$. But weights cannot prove it. The six sets of Figure~\ref{fig:tiling15} cover every point of $[15]$ exactly twice, so giving each of them the value $\frac12$ is a fractional tiling with total $3$, and by~\eqref{eq:tiling} every weighting that gives each $5$-element Sidon set weight at least $L$ has total weight at least $3L$. A proof that $\SR(3)\leq15$ must use that a class is taken whole or not at all; an exhaustive search does.
\end{example}

\subsection{Nine colors}\label{sec:nine}

Let
\begin{align*}
 A&=\{1,8,9,18,22,37,48,64,70,82,102,105,107\},\\
 B&=\{1,3,6,26,38,44,60,71,86,90,99,100,107\},
\end{align*}
and put $A_i=A+i$ and $B_i=B+i$ for $i=0,1,2$; these six sets lie in $[109]$. The reflection $\rho(x)=110-x$ maps $[109]$ to itself, and one checks from the lists that
\begin{equation}\label{eq:reflection}
 \rho(A_i)=B_{2-i},\qquad \rho(B_i)=A_{2-i}.
\end{equation}
For $C\subset[109]$, let $\SSS(C)$ be the family of all $12$-element Sidon subsets of $[109]\setminus C$.

\begin{table}[ht]
\centering
\begin{tabular}{lrrrr}
\toprule
Removed set $C$&$|\SSS(C)|$&Weights&Total weight&Smallest weight of a set\\
\midrule
$A_0$&191515&$w_0$&79544&10003\\
$A_1$&185834&$w_1$&79565&10003\\
$A_2$&174522&$w_2$&79663&10003\\
$A_0\cup B_1$&7285&$w'$&52469&10002\\
\bottomrule
\end{tabular}
\caption{The four families of $12$-element Sidon sets used in Step~3 of the proof, and their weight checks.}
\label{tab:counts}
\end{table}

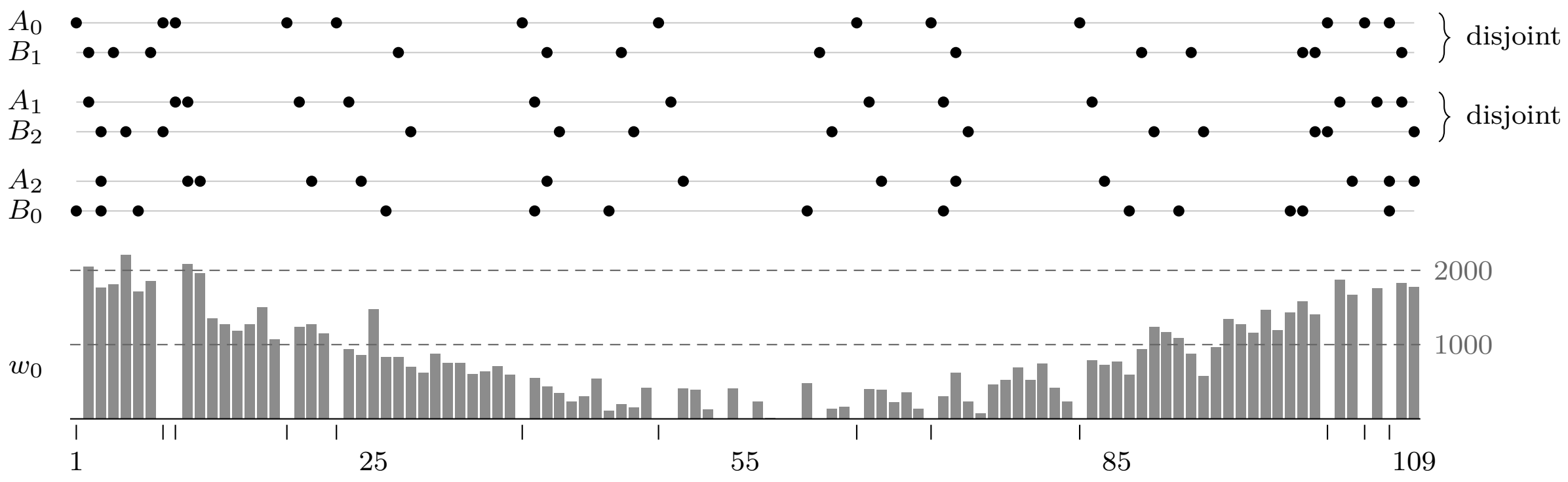
\begin{figure}[ht]
\centering
\begin{tikzpicture}[x=0.125cm,y=0.5cm]
\node[font=\scriptsize,anchor=east] at (-0.5,0) {$A_0$};
\draw[black!20] (1,0) -- (109,0);
\foreach \x in {1,8,9,18,22,37,48,64,70,82,102,105,107} \fill (\x,0) circle (1.6pt);
\node[font=\scriptsize,anchor=east] at (-0.5,-0.6) {$B_1$};
\draw[black!20] (1,-0.6) -- (109,-0.6);
\foreach \x in {2,4,7,27,39,45,61,72,87,91,100,101,108} \fill (\x,-0.6) circle (1.6pt);
\node[font=\scriptsize,anchor=east] at (-0.5,-1.6) {$A_1$};
\draw[black!20] (1,-1.6) -- (109,-1.6);
\foreach \x in {2,9,10,19,23,38,49,65,71,83,103,106,108} \fill (\x,-1.6) circle (1.6pt);
\node[font=\scriptsize,anchor=east] at (-0.5,-2.2) {$B_2$};
\draw[black!20] (1,-2.2) -- (109,-2.2);
\foreach \x in {3,5,8,28,40,46,62,73,88,92,101,102,109} \fill (\x,-2.2) circle (1.6pt);
\node[font=\scriptsize,anchor=east] at (-0.5,-3.2) {$A_2$};
\draw[black!20] (1,-3.2) -- (109,-3.2);
\foreach \x in {3,10,11,20,24,39,50,66,72,84,104,107,109} \fill (\x,-3.2) circle (1.6pt);
\node[font=\scriptsize,anchor=east] at (-0.5,-3.8) {$B_0$};
\draw[black!20] (1,-3.8) -- (109,-3.8);
\foreach \x in {1,3,6,26,38,44,60,71,86,90,99,100,107} \fill (\x,-3.8) circle (1.6pt);
\draw[decorate,decoration={brace,amplitude=3pt}] (111,0.2) -- (111,-0.8) node[midway,right=4pt,font=\scriptsize] {disjoint};
\draw[decorate,decoration={brace,amplitude=3pt}] (111,-1.4) -- (111,-2.4) node[midway,right=4pt,font=\scriptsize] {disjoint};
\node[font=\scriptsize,anchor=east] at (-0.5,-7.00) {$w_0$};
\foreach \x/\h in {2/3.078,3/2.654,4/2.718,5/3.311,6/2.567,7/2.790,10/3.134,11/2.943,12/2.037,13/1.915,14/1.784,15/1.911,16/2.256,17/1.605,19/1.856,20/1.914,21/1.730,23/1.407,24/1.286,25/2.216,26/1.248,27/1.244,28/1.046,29/0.938,30/1.320,31/1.134,32/1.134,33/0.903,34/0.963,35/1.068,36/0.896,38/0.825,39/0.652,40/0.519,41/0.345,42/0.458,43/0.808,44/0.168,45/0.297,46/0.233,47/0.624,50/0.618,51/0.588,52/0.194,54/0.621,56/0.355,57/0.015,59/0.006,60/0.717,62/0.203,63/0.247,65/0.603,66/0.591,67/0.333,68/0.537,69/0.200,71/0.450,72/0.938,73/0.350,74/0.108,75/0.696,76/0.790,77/1.044,78/0.784,79/1.113,80/0.630,81/0.355,83/1.190,84/1.091,85/1.157,86/0.891,87/1.409,88/1.857,89/1.748,90/1.639,91/1.311,92/0.873,93/1.449,94/2.021,95/1.913,96/1.744,97/2.198,98/1.794,99/2.151,100/2.377,101/2.107,103/2.814,104/2.513,106/2.635,108/2.751,109/2.670} \fill[black!45] (\x-0.4,-8.0) rectangle (\x+0.4,-8.0+\h);
\draw (0.5,-8.0) -- (109.5,-8.0);
\foreach \x in {1,8,9,18,22,37,48,64,70,82,102,105,107} \draw (\x,-8.12) -- (\x,-8.42);
\draw[densely dashed,black!60] (0.5,-5.00) -- (109.5,-5.00) node[right,font=\scriptsize] {$2000$};
\draw[densely dashed,black!60] (0.5,-6.50) -- (109.5,-6.50) node[right,font=\scriptsize] {$1000$};
\node[font=\scriptsize,below] at (1,-8.40) {$1$};
\node[font=\scriptsize,below] at (25,-8.40) {$25$};
\node[font=\scriptsize,below] at (55,-8.40) {$55$};
\node[font=\scriptsize,below] at (85,-8.40) {$85$};
\node[font=\scriptsize,below] at (109,-8.40) {$109$};
\end{tikzpicture}
\caption{Top: the six $13$-element Sidon subsets of $[109]$, one per row. The sets $A_i$ are translates of one ruler and the sets $B_i$ of its mirror image; the only disjoint pairs are $\{A_0,B_1\}$ and $\{A_1,B_2\}$. Bottom: the weights $w_0$ of Step~3, which exclude eight disjoint $12$-element Sidon sets outside $A_0$ (the ticks below the axis mark the points of $A_0$, which are removed). As in Example~\ref{ex:weights8}, the points near the ends are the expensive ones.}
\label{fig:nine}
\end{figure}

\begin{proof}[Proof that $\SR(9)\leq109$]
Suppose that $[109]$ is partitioned into nine Sidon sets.

\emph{Step 1: the classes with $13$ elements.} An exhaustive enumeration (Section~\ref{sec:computation}) shows that every Sidon subset of $[109]$ has at most $13$ elements, that the $13$-element Sidon subsets are exactly $A_0,A_1,A_2,B_0,B_1,B_2$, and that among them the only disjoint pairs are $\{A_0,B_1\}$ and $\{A_1,B_2\}$ (Figure~\ref{fig:nine}, top); in particular, no three are pairwise disjoint. The slack is $9\cdot13-109=8$ (Figure~\ref{fig:slack}), so at least one class of the partition has $13$ elements, and at most two do.

\emph{Step 2: the other classes.} If exactly one class has $13$ elements, it is some $A_i$ or $B_i$, and the other eight classes have at most $12$ elements and partition its complement, which has $96$ points. Their slack is $8\cdot12-96=0$, so all eight have exactly $12$ elements, and they are eight pairwise disjoint members of $\SSS(A_i)$ or $\SSS(B_i)$. If two classes have $13$ elements, they are $A_0,B_1$ or $A_1,B_2$, and the other seven classes partition the remaining $83$ points. Their slack is $7\cdot12-83=1$, so at least six of them have exactly $12$ elements, and they are six pairwise disjoint members of $\SSS(A_0\cup B_1)$ or $\SSS(A_1\cup B_2)$.

\emph{Step 3: the weights.} For each row of Table~\ref{tab:counts}, the column $w_0$, $w_1$, $w_2$ or $w'$ of Table~\ref{tab:weights} (Appendix~\ref{app:weights}) gives every member of $\SSS(C)$ weight at least $10000$: the last column of Table~\ref{tab:counts} is the smallest weight over the complete list of $\SSS(C)$, and the total weight is the column sum. Since $8\cdot10000>79663$ and $6\cdot10000>52469$, Lemma~\ref{lem:weights}\ref{it:pack} with $L=10000$ shows that for $i=0,1,2$ at most seven members of $\SSS(A_i)$ are pairwise disjoint, and at most five members of $\SSS(A_0\cup B_1)$ are pairwise disjoint. Because of~\eqref{eq:reflection}, the same holds for $\SSS(B_0)$, $\SSS(B_1)$, $\SSS(B_2)$ and $\SSS(A_1\cup B_2)$: if $w$ works for $C$, then $x\mapsto w(110-x)$ works for $\rho(C)$, since the reflection preserves Sidon sets (Section~\ref{sec:maps}), disjointness and total weight. This contradicts Step~2.

Every finite fact used here is checked by the program described in Section~\ref{sec:computation}.
\end{proof}

The weights look like those of Example~\ref{ex:weights8} (Figure~\ref{fig:nine}, bottom): heavy near both ends of $[109]$ and light in the middle. The reason is the same. Since $g_{12}=85$, a $12$-element Sidon subset of $[109]$ has its smallest element at most $24$ and its largest at least $86$, so every candidate must pay for points near both ends.

\subsection{Ten colors}\label{sec:ten}

For ten colors the plan needs one more round: after the largest classes are fixed, the slack of the remaining classes is $0$ or $1$, so what is left must be split almost exactly into $13$-element Sidon sets. Lemma~\ref{lem:weights}\ref{it:cover}, with weights of both signs, excludes this.

\begin{proof}[Proof that $\SR(10)\leq132$]
Since $g_{15}=151>131$, every Sidon subset of $[132]$ has at most $14$ elements, and the slack of a Sidon $10$-partition of $[132]$ would be $10\cdot14-132=8$. So it has $t\geq2$ classes with $14$ elements. These are among the $52$ translates of the $24$ rulers of Table~\ref{tab:near} with $14$ marks; $134$ pairs and $8$ triples of them are pairwise disjoint, and no four, so $t=2$ or $t=3$. Let $U$ be the set of points outside the $t$ classes with $14$ elements. The other $10-t$ classes have at most $13$ elements and cover $U$.
\begin{itemize}
\item If $t=2$, then $|U|=104$ and the slack of the eight remaining classes is $8\cdot13-104=0$. So all of them have $13$ elements, and $U$ is a disjoint union of $13$-element Sidon sets.
\item If $t=3$, then $|U|=90$ and the slack of the seven remaining classes is $7\cdot13-90=1$. So six of them have $13$ elements and one has $12$, and $U$ is the disjoint union of six $13$-element Sidon sets and $12$ further points.
\end{itemize}
For each of the $142$ possible sets $U$, let $\mathcal F$ be the family of $13$-element Sidon subsets of $[132]$ contained in $U$; it has between $14{,}116$ and $162{,}303$ members when $t=2$ and between $268$ and $2{,}363$ when $t=3$. For each $U$ the verification repository~\cite{SRVerif} gives integers $w(x)$, with $|w(x)|\leq1000$, that satisfy Lemma~\ref{lem:weights}\ref{it:cover} with $f=0$ if $t=2$ and $f=12$ if $t=3$. So no Sidon $10$-partition of $[132]$ exists.
\end{proof}

\subsection{Why the method stops at \texorpdfstring{$109$}{109}}\label{sec:stop}

For $[108]$ the argument of Section~\ref{sec:nine} cannot work, for the reason seen in Example~\ref{ex:tiling15}. Since $9\cdot12=108$, a Sidon $9$-partition of $[108]$ could consist of nine classes with $12$ elements each (slack $0$ with $M=12$), so one has to exclude nine pairwise disjoint $12$-element Sidon sets in $[108]$. There are $1{,}213{,}118$ such sets, and $108$ of them carry nonnegative rational numbers $y_S$ with $\sum_{S\ni x}y_S=1$ for every $x\in[108]$, so that $\sum_S y_S=9$. This fractional tiling is in the verification repository, with a program that checks it in exact arithmetic. By~\eqref{eq:tiling}, any weights that give every $12$-element Sidon subset of $[108]$ weight at least $L$ have total weight at least $9L$, so Lemma~\ref{lem:weights}\ref{it:pack} can never exclude nine pairwise disjoint members: the linear program of Section~\ref{sec:weights} has value exactly $9$ for this family. Going below $109$ requires an argument that uses integrality, such as an exact search.

\subsection{The computations}\label{sec:computation}

\paragraph{Enumerating Sidon sets.} To list the $r$-element Sidon subsets of a set $X\subset[n]$, we build each set in increasing order. Along with the current set $S=\{a_1<\cdots<a_s\}$ we keep its set of differences $a_j-a_i$. A candidate $x>a_s$ adds the differences $x-a_1,\ldots,x-a_s$, which are different from one another, so the set $S\cup\{x\}$ is Sidon exactly when none of them is already present. The search tries every admissible $x$ in increasing order, outputs a set when it reaches $r$ elements, and abandons a branch when too few candidates remain. Every Sidon set is found exactly once, along the increasing sequence of its elements, because every prefix of a Sidon set is Sidon. For the whole interval it is enough to list the sets that contain $1$ and translate them. This first implementation also uses lower bounds for ruler lengths to cut branches (Section~\ref{sec:upper-notes}).

A second implementation keeps the pair sums instead of the differences (including $2a$ for $a\in S$) and uses no ruler bounds. Both programs return identical lists, compared as sorted lists of bit masks, not only identical counts. For Step~1 of the proof for nine colors we list the subsets of $[109]$ with $13$ and $14$ elements, and for Step~3 the $12$-element subsets of the four complements in Table~\ref{tab:counts}; smaller sets are never listed.

\paragraph{Finding the weights.} For a family $\mathcal F=\SSS(C)$ we solve the linear program
\begin{equation}\label{eq:lp}
 \text{minimize }\sum_{x\in[109]}y_x\quad\text{subject to}\quad
 \sum_{x\in S}y_x\geq1\ (S\in\mathcal F),\qquad y_x\geq0,
\end{equation}
starting with a subfamily and adding violated constraints until none remain, with an off-the-shelf solver (Section~\ref{sec:upper-notes}). We multiply the solution by $10000$, round up, and check the resulting integers against every set of the family. The signed weights for ten colors come from the corresponding linear program for Lemma~\ref{lem:weights}\ref{it:cover}, solved and rounded in the same way. Only the integers enter the proofs.

\paragraph{Checking the certificates.} A program in the verification repository~\cite{SRVerif} regenerates the lists of $12$-, $13$- and $14$-element Sidon sets used in the proof for nine colors and checks Step~1 (the six $13$-element sets, the absence of $14$-element ones, and the disjoint pairs), the reflection identities~\eqref{eq:reflection} and, for each row of Table~\ref{tab:counts}, the number of sets, that every set has weight at least $10000$, and the total weight. The verification has not been formalized in a proof assistant. Appendix~\ref{app:computations} lists the commands and their running times.

\paragraph{Rulers, and ten to twelve colors.} The rulers of Table~\ref{tab:near}, and the $3{,}423{,}192$ rulers with $13$ marks and length at most $131$ from which the families $\mathcal F$ of Section~\ref{sec:ten} are formed, were listed by a third program, independent of the two above. It places the marks in increasing order, keeping the set of differences used so far, and abandons a branch when the remaining marks cannot fit, either because a ruler with that many marks is longer than the space left or because the smallest unused differences already add up to more. The first test uses lower bounds for $g_r$, and the program proves each one before using it: running it for $r=3,\ldots,16$ marks in turn recovers $g_3,\ldots,g_{16}$, the shortest length found in each run, and it also finds no ruler with $17$ marks and length at most $184$. Its lists agree with those of two earlier implementations, and it reproduces the numbers of Sidon sets in Sections~\ref{sec:computation} and~\ref{sec:exact} (for example $1{,}213{,}118$ twelve-element Sidon subsets of $[108]$ and $125{,}628{,}796$ eleven-element ones of $[102]$).

For ten colors, the condition of Lemma~\ref{lem:weights}\ref{it:cover} involves the sum of the $f$ smallest weights. Instead of sorting, each certificate stores an integer $z$ besides $w$, and the program checks that
\[
 w(U)<fz-\sum_{x\in U}\max\{0,z-w(x)\}.
\]
This implies the condition of the lemma: every point satisfies $w(x)\geq z-\max\{0,z-w(x)\}$, so any $f$ points of $U$ have total weight at least the right-hand side. (Taking for $z$ the $f$-th smallest weight, the right-hand side is exactly the sum of the $f$ smallest weights.) The same program lists the rulers, checks the disjoint families of Section~\ref{sec:upper1112} and checks every weight against the complete families $\mathcal F$. As a second check, which does not use the weights, an exhaustive exact cover search over the same families finds no decomposition of any of the $142$ sets $U$; with the four $12$-element classes of our partition of $[102]$ fixed, the same program finds exactly that partition.

\subsection{Notes}\label{sec:upper-notes}

The bound $13$ in Step~1 for nine colors and the description of the $13$-element Sidon subsets of $[109]$ by translates agree with Espinosa-Garc\'ia and Pellicer~\cite{EP}; the bound also follows from $g_{14}=127$, see~\cite[Table~1]{EP}. The first enumeration program cuts branches with the lower bounds for ruler lengths in~\cite[Table~1]{EP}. The linear programs~\eqref{eq:lp} were solved with GLOP from OR-Tools~\cite{ORTools}. Lemma~\ref{lem:weights}\ref{it:pack} is the easy half of linear programming duality, which in this setting relates fractional packings and fractional covers of a hypergraph.

\section{Lower bounds by computer search}\label{sec:search}

A lower bound $\SR(k)>N$ is a certificate: once a Sidon $k$-partition of $[N]$ is found, checking it is easy, and how it was found does not matter for its correctness. This section first proves the lower bounds of Theorem~\ref{thm:main} by checking the partitions, and then explains how they were found.

\begin{proof}[Proof of the lower bounds in Theorem~\ref{thm:main}]
For each $k$, the verification repository~\cite{SRVerif} contains a Sidon $k$-partition of $[N]$, where $N+1$ is the current lower bound in Table~\ref{tab:bounds}. To check that a class $A$ is a Sidon set it suffices to compute its $\binom{|A|+1}2$ sums $a+b$ with $a\leq b$ and see that they are all different. A program in the repository does this for every class of every partition, checks that the classes partition $[N]$, and, as a second test, checks that the differences in each class are distinct (Appendix~\ref{app:computations}).
\end{proof}

\subsection{Counting conflicts}

A search needs a measure of how far a coloring is from being a Sidon partition. For a coloring $c:[N]\to[k]$ and a color $j$, let $\nu_j(\sigma)$ be the number of pairs $\{a,b\}$ with $a\leq b$, $a+b=\sigma$ and $c(a)=c(b)=j$. The number of \emph{conflicts} is
\begin{equation}\label{eq:conflicts}
 C(c)=\sum_{j=1}^k\sum_{\sigma=2}^{2N}\binom{\nu_j(\sigma)}2,
\end{equation}
the number of equalities $a+b=a'+b'$ between two different pairs of the same color. A coloring is a Sidon partition exactly when $C(c)=0$. For example, putting $8$ in the white class of Figure~\ref{fig:sr2} gives the coloring $\{1,2,5,7\}\cup\{3,4,6,8\}$ of $[8]$, which has one conflict, $4+8=6+6$.

A search changes one point at a time, so it needs the change in $C$, not $C$ itself. When a point $x$ moves from class $A_i$ to class $A_j$, the pairs that leave $A_i$ have the sums $x+a$ ($a\in A_i$, including $a=x$) and the pairs that enter $A_j$ have the sums $x+b$ ($b\in A_j$) and $2x$. These sums are different from one another, so, with the counts taken before the move,
\begin{equation}\label{eq:delta}
 \Delta C=-\sum_{a\in A_i}\bigl(\nu_i(x+a)-1\bigr)
             +\nu_j(2x)+\sum_{b\in A_j}\nu_j(x+b).
\end{equation}
This takes $O(N/k)$ operations, instead of the $O(N^2)$ needed to recount.

\subsection{Extension and repair}\label{sec:repair}

\begin{figure}[ht]
\centering
\begin{tikzpicture}[x=0.82cm]
\node[font=\scriptsize,anchor=east] at (0.4,0) {before};
\foreach \x/\s in {1/clA,2/clA,3/clB,4/clC,5/clA,6/clC,7/clA,8/clB,9/clB,10/clC,11/clC,12/clB} \node[\s] (a\x) at (\x,0) {$\x$};
\node[clA,dashed,draw=black] (a13) at (13,0) {$13$};
\draw[thick] (a1.north) to[bend left=22] node[above,font=\scriptsize] {$6$} (a7.north);
\draw[thick] (a7.north) to[bend left=22] node[above,font=\scriptsize] {$6$} (a13.north);
\node[font=\scriptsize,anchor=east] at (0.4,-1.4) {after};
\foreach \x/\s in {1/clB,2/clA,3/clB,4/clC,5/clA,6/clC,7/clA,8/clB,9/clB,10/clC,11/clC,12/clB,13/clA} \node[\s] (b\x) at (\x,-1.4) {$\x$};
\draw[->,densely dashed] (a1.south) -- (b1.north);
\end{tikzpicture}
\caption{Extension and repair with three colors. Before: the Sidon $3$-partition $\{1,2,5,7\}\cup\{3,8,9,12\}\cup\{4,6,10,11\}$ of $[12]$ (dark, gray, white), extended by $13$. The point $13$ creates a conflict in every class: $1+13=7+7$ in the dark class, $3+13=8+8$ and $8+13=9+12$ in the gray one, and $4+13=6+11$ in the white one. It is placed in the dark class, where it creates one conflict (the repeated difference $6$). After: moving $1$ to the gray class removes the conflict, and creates none, so $\{2,5,7,13\}\cup\{1,3,8,9,12\}\cup\{4,6,10,11\}$ is a Sidon $3$-partition of $[13]$.}
\label{fig:repair}
\end{figure}
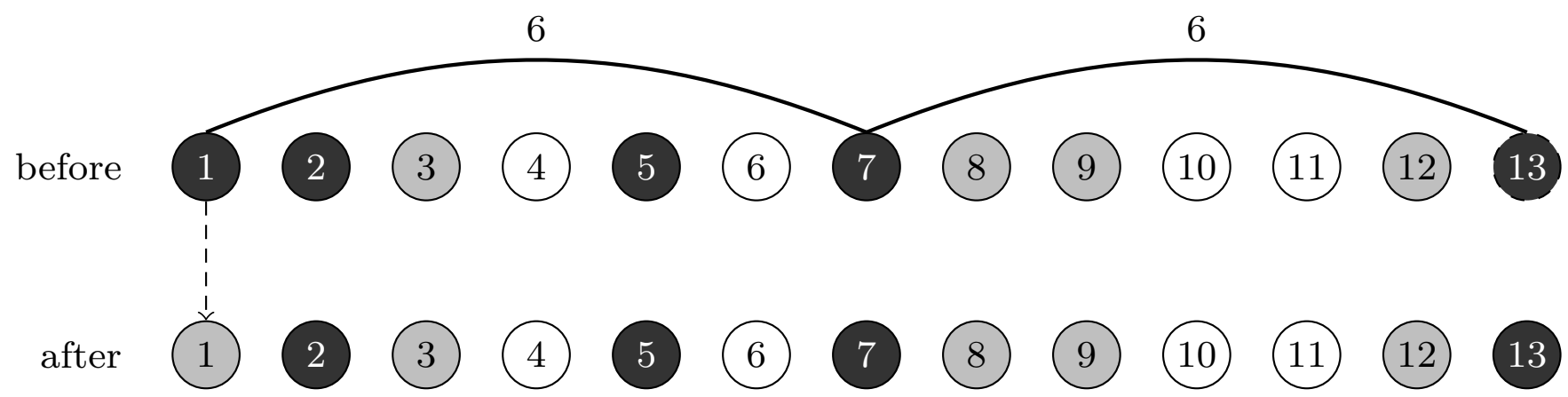

The idea is to grow a partition one point at a time. A Sidon partition of $[N-1]$ is almost a partition of $[N]$: only the new point is missing, and wherever it is placed it creates a few conflicts. A local search then tries to remove them by recoloring a few points, without disturbing the rest. Figure~\ref{fig:repair} shows the smallest kind of case, where a single move repairs the partition.

A search starts from a Sidon partition of $[N-1]$, or of a shorter interval, and places each new point $x$ in a class $A_j$ that creates the fewest conflicts, $\nu_j(2x)+\sum_{y\in A_j}\nu_j(x+y)$, breaking ties at random. Then it repairs. Each step picks a random conflict, considers moving each of its points to every other class, and makes the move with the smallest $\Delta C$ among those allowed. Left alone, such a search tends to cycle, undoing its own moves, or to get stuck where every move makes things worse. Three standard devices prevent this:
\begin{itemize}
\item a point may not move back to a class it left in the last few steps (a \emph{tabu list});
\item a move that increases $C$ may be accepted, with a probability that grows the longer the search goes without improving, as in simulated annealing;
\item in some runs, conflicts that survive for a long time get larger weights, so that the search concentrates on them.
\end{itemize}
The search also uses two moves that look further ahead. A \emph{preparatory move} recolors a point without changing the set of conflicts, but so that a later move can remove one of them more cheaply. A \emph{chain} is a short sequence of moves that creates new conflicts on the way and removes them by the end. When nothing improves for a long time, a \emph{mutation} recolors a few random points, and the repair starts again. Each partition found becomes the start for the next interval.

\subsection{Starting colorings}

For eight colors we started from a published partition of $[80]$ (Section~\ref{sec:search-notes}); those searches stopped at two conflicts on $[81]$, and with nine colors no search reached $[100]$, which is why we turned to the exact search of Section~\ref{sec:exact}.

For the local search with nine colors we started from nine published disjoint $11$-element Sidon sets covering $[100]\setminus\{46\}$, and for ten colors from a published partition of $[110]$ into $11$-element sets and from ten published disjoint $12$-element Sidon sets covering $[123]\setminus\{44,73,92\}$; we restricted these to the target interval and colored the missing points. For eleven colors we started from our partition of $[118]$ into ten classes, with an empty eleventh class. For $14$ and more colors we started from the cyclic partitions of Section~\ref{sec:cyclic}. For $17$, $18$ and $20$ colors we started thousands of searches at once, each from a different affine image of the Singer partition, on a GPU.

\subsection{Deleting classes}\label{sec:delete}

For $9$, $13$, $15$ and $16$ colors the best partitions of the local search came from cyclic partitions with too many classes. Take an affine image of a Singer partition with $q+1>k$ classes, delete $q+1-k$ of its classes, and restrict what remains to an interval $[N]$ well below $v$. The remaining classes are Sidon sets, so this start has no conflicts at all; what is missing are the points of the deleted classes, which start uncolored. Thousands of such starts were placed and repaired at once on a GPU, extending one point at a time as in Section~\ref{sec:repair}, and the best partitions reached there were extended further by the CPU search:
\begin{center}
\begin{tabular}{rrrrr}
\toprule
$k$&Singer $q$&classes deleted&GPU&CPU\\
\midrule
9&13&5&$[96]$&$[99]$\\
13&13&1&$[170]$&$[174]$\\
15&16&2&$[217]$&$[223]$\\
16&16&1&$[244]$&$[250]$\\
\bottomrule
\end{tabular}
\end{center}
By contrast, extending the previous partition of $[243]$ with sixteen colors had stayed at one conflict on $[244]$ for six hours, and every earlier search on $[99]$ with nine colors had stopped at two conflicts; from the deleted-class start, the CPU search went from $[244]$ to $[250]$ in an hour and a half. Deleting one class from our best partition with $k+1$ colors (read on its subintervals and their reflections, which preserve Sidon sets, Section~\ref{sec:maps}) also works: it reached the same intervals for $9$, $11$, $13$ and $15$ colors, in minutes, but not beyond them.

\subsection{Exact search with the largest classes fixed}\label{sec:exact}

The partitions for eight and nine colors were found by our exact search. It turns the rigidity of Section~\ref{sec:upper} into a search strategy: a partition must contain several classes of the largest size, and since there are few such sets, we choose them first. What remains is a much smaller and more constrained problem, which we solve exhaustively. We describe the search for $[81]$ with eight colors; for $[102]$ with nine colors every size is one larger. Its programs are in the verification repository~\cite{SRVerif}, with the commands of the runs described below and positive controls.

\paragraph{Which classes to fix.} A Sidon subset of $[81]$ has at most $11$ elements. The slack of a Sidon $8$-partition of $[81]$ is $8\cdot11-81=7$, so at least $t\geq1$ classes have $11$ elements. We split the problem by $t$ and treat each case as a choice of $t$ disjoint $11$-element Sidon sets for the large classes. The complement of these sets must then be split into Sidon sets with at most $10$ elements. The more large classes are fixed, the smaller and more constrained this remaining problem is, so we started with the largest values of~$t$.

\paragraph{The search.} For a fixed choice of large classes we do not color the remaining points one at a time; we choose whole classes.
\begin{enumerate}
\item List all Sidon sets of the two largest sizes, as in Section~\ref{sec:computation}. $[81]$ has $4030$ Sidon sets with $11$ elements and $11{,}817{,}176$ with $10$ elements.
\item List the $t$-tuples of pairwise disjoint $11$-element Sidon sets, up to the reflection $x\mapsto82-x$. For $t=3$ there are $621{,}365$ of them.
\item For each $t$-tuple, keep only the $10$-element Sidon sets that avoid it. The remaining $81-11t$ points must be covered by some of these sets together with a few smaller classes, whose sizes are fixed in advance. For $t=3$ there are five more classes on $48$ points; their slack is $5\cdot10-48=2$, so their sizes are $10,10,10,10,8$ or $10,10,10,9,9$, and we ran both. The smaller classes are not chosen from a list; their points are collected in one \emph{leftover} set.
\item Cover the points exactly. At each step take the uncovered point with the fewest options and try each of them. An option is a remaining $10$-element set that contains the point, or putting the point in the leftover. The leftover option is allowed only while the leftover can still be split into Sidon sets of the prescribed sizes.
\item When all the $10$-element classes are chosen, a small exact search splits the leftover into Sidon sets of the prescribed sizes.
\end{enumerate}
For the triples we sampled, the complement of three disjoint $11$-element Sidon sets contains about $1700$ of the $11{,}817{,}176$ $10$-element Sidon sets (the median), and a triple takes about $0.1$ seconds on average. Allowing $9$- and $8$-element sets as candidates as well gives millions of candidates and was much slower than pooling them in the leftover.

\paragraph{Nine colors.} A Sidon subset of $[102]$ has at most $12$ elements, and the slack $9\cdot12-102=6$ gives $t\geq3$ classes with $12$ elements. With $t=4$, the other five classes hold $54$ points with slack $5\cdot11-54=1$, so their sizes are $11,11,11,11,10$, and the candidates are the $125{,}628{,}796$ Sidon sets in $[102]$ with $11$ elements. The same search on $[100]$ and $[101]$ found the partitions that first gave $\SR(9)\geq101$ and $\SR(9)\geq102$.

\paragraph{Results.} Table~\ref{tab:exact} summarizes the runs.
\begin{table}[ht]
\centering
\begin{tabular}{llrrl}
\toprule
Interval&Large classes fixed&Cases&Cases run&Result\\
\midrule
$[81]$, $k=8$&three of size $11$&$621{,}365$&all&exactly $5$ partitions\\
$[100]$, $k=9$&four of size $12$&$4{,}622{,}302$&$246{,}575$&first partition\\
$[101]$, $k=9$&four of size $12$&$86{,}166{,}470$&about $1.29$ million&first partition\\
$[102]$, $k=9$&four of size $12$&about $1.5\times10^9$&about $0.9\times10^9$&first partition\\
\bottomrule
\end{tabular}
\caption{The exact search. A case is a $t$-tuple of disjoint Sidon sets of the largest size, up to reflection. For $[100]$, $[101]$ and $[102]$, ``Cases run'' is the number of cases run before the first partition appeared.}
\label{tab:exact}
\end{table}
For $[81]$ the run is exhaustive: no case timed out, and it found exactly five triples that can be completed, each in exactly one way and with class sizes $11,11,11,10,10,10,9,9$. Up to reflection these are all the Sidon $8$-partitions of $[81]$ with exactly three $11$-element classes, and no two of them share a class, even after reflecting one of them. A separate exhaustive search over the $12{,}918$ ($6{,}483$ up to reflection) disjoint quadruples of $11$-element Sidon sets found no partition with four or more such classes. Partitions with one or two $11$-element classes were not examined, because each case has far more candidates.

The partition of $[101]$ has class sizes $12,12,12,12,11,11,11,10,10$. The number of $4$-tuples of $[102]$ in Table~\ref{tab:exact} was estimated from a sample. On $[102]$ all the $5$- and $6$-tuples were also run, with every size profile in which the other classes have at most $11$ elements, and none can be completed. Moreover, $[102]$ has no seven pairwise disjoint $12$-element Sidon sets, by a separate exhaustive search that shares no code with the exact cover. Hence no Sidon $9$-partition of $[102]$ has five or more classes with $12$ elements. In hindsight, the partition of $[101]$ already gives $\SR(9)\geq103$: adding $102$ to its class $\{2,10,17,48,49,51,65,74,78,84\}$ keeps that class a Sidon set, and the run on $[102]$ found this extended partition too. The partition of $[102]$ in the verification repository is the other one; it shares no class with the partition of $[101]$.

\paragraph{Tests.} We tested the search on cases with a known answer. On $[80]$, with some classes of known partitions fixed, it finds those partitions again, and its solution counts agree with those of an independent program that completes partial colorings point by point. With the sizes lowered by one, the nine-color program recovers our partition of $[81]$ from its three $11$-element classes in under a second.

\subsection{Open targets}

Table~\ref{tab:open} shows, for each number of colors, the first interval for which we have no partition and the fewest conflicts the searches reached there. Typically many different colorings with a single conflict were found, and for every $k$ from $10$ to $19$ except $12$ the extensions that reached the bound in Table~\ref{tab:bounds} stopped at a single conflict on the next interval; for $12$ colors they stopped at two.

\begin{table}[ht]
\centering
\begin{tabular}{rrr@{\qquad}rrr@{\qquad}rrr}
\toprule
$k$&Target&Conflicts&$k$&Target&Conflicts&$k$&Target&Conflicts\\
\midrule
8&$[82]$&--&12&$[155]$&2&16&$[251]$&1\\
9&$[103]$&--&13&$[175]$&1&17&$[288]$&1\\
10&$[119]$&1&14&$[196]$&1&18&$[324]$&1\\
11&$[134]$&1&15&$[224]$&1&19&$[352]$&1\\
\bottomrule
\end{tabular}
\caption{The first interval without a known Sidon $k$-partition, and the smallest number of conflicts~\eqref{eq:conflicts} reached there (``--'': not recorded).}
\label{tab:open}
\end{table}

Appendix~\ref{app:computations} records the machines the searches ran on, and other methods that were tried on these targets without success.

\subsection{Notes}\label{sec:search-notes}

The extension-and-repair approach is adapted from Cremer~\cite{Cremer}, who improved the lower bounds for the Ramsey numbers $R(6,8)$ and $R(8,10)$ by taking known good colorings, adding vertices, and removing the resulting monochromatic cliques by local search. Our implementation is our own; the preparatory moves and chains follow his. Repairing a random conflict with the best move is the min-conflicts heuristic of Minton et al.~\cite{Minton}; the three devices are tabu search~\cite{HertzDeWerra}, simulated annealing~\cite{Kirkpatrick} and the breakout method~\cite{Morris}. The exact cover step of Section~\ref{sec:exact} branches like Knuth's Algorithm~X~\cite{Knuth}. Fixing most of a coloring and solving exactly for the rest, as in the unsuccessful searches of Appendix~\ref{app:computations}, is large neighborhood search~\cite{Shaw}; we used PySAT~\cite{PySAT} with the SAT solver CaDiCaL~\cite{CaDiCaL}. The published starting colorings are those of Xiu, Fan and Liang~\cite{XFL}: for eight colors the partition of $[80]$ in their Table~I, and for nine and ten colors the sets in their Table~IV.

\section{Van der Waerden numbers}\label{sec:vdw}

In this section a \emph{progression} is a $3$-term arithmetic progression $x,x+d,x+2d$ with $d\geq1$, and an $r$-coloring of $[n]$ is \emph{valid} if no progression in $[n]$ is monochromatic. So $W(r,3)>n$ exactly when $[n]$ has a valid $r$-coloring. Progression-free sets are much denser than Sidon sets, so valid colorings are long (thousands of points for $13$ to $17$ colors), and exhaustive methods are out of reach. Our colorings come from constructions that turn short colorings into long ones (Sections~\ref{sec:vdw-derived} and~\ref{sec:vdw1516}), improved by a local search (Section~\ref{sec:vdw-local}). Table~\ref{tab:vdw-published} lists the published lower bounds that we use.

\begin{table}[ht]
\centering
\begin{tabular}{rrl@{\qquad}rrl}
\toprule
$r$&$W(r,3)>$&Source&$r$&$W(r,3)>$&Source\\
\midrule
5&170&Heule~\cite{Heule}&9&763&Heule~\cite{Heule}\\
6&225&BCT~\cite{BCT}&10&892&Komkov~\cite{Komkov}\\
7&343&Komkov~\cite{Komkov}&11&1187&Komkov~\cite{Komkov}\\
8&515&Komkov~\cite{Komkov}&12&2329&Heule~\cite{Heule}\\
&&&17&3549&Komkov~\cite{Komkov}\\
\bottomrule
\end{tabular}
\caption{Published lower bounds for $W(r,3)$, $5\leq r\leq17$. The value for six colors comes from the recurrence of~\cite{BCT}; Heule's coloring gives $223$. The table of~\cite{BCT} lists $W(9,3)>775$, without a source; we use Heule's $763$, which does not affect anything below. For $17$ colors, Heule's repository~\cite{HeuleRepo} has a coloring showing $W(17,3)>3546$, which Komkov raised to $3549$.}
\label{tab:vdw-published}
\end{table}

\subsection{Combining colorings}\label{sec:vdw-derived}

There are two ways to build a valid coloring from smaller ones. The first puts two valid colorings side by side, with disjoint sets of colors: a progression in $[n+n']$ that meets both blocks has two colors, so
\begin{equation}\label{eq:superadditive}
 W(r+s,3)-1\geq\bigl(W(r,3)-1\bigr)+\bigl(W(s,3)-1\bigr).
\end{equation}

The second is a product, and it gains much more. Take a valid coloring of $\{0,\ldots,n-1\}$ with only $s$ colors and blow up every point $t$ into a block of $q$ consecutive points $qt,qt+1,\ldots,qt+q-1$. The point $qt+u$ of the block, with $0\leq u<q$, gets its color from a \emph{palette} of $s$ colors attached to the residue $u$; which entry of the palette it gets is decided by the color of $t$ in the short coloring. A progression in the long coloring either changes residue as it goes, or stays on one residue. In the first case it is a progression modulo $q$, and we choose the palettes so that every color appears only on a set of residues without progressions modulo $q$. In the second case it is a blown-up copy of a progression of the short coloring, which was not monochromatic.

To make this precise, write $r_3(\Z_q)$ for the largest size of a subset of $\Z_q$ that contains no progression $x,x+d,x+2d$ modulo $q$ with $d\not\equiv0$; when $q$ is even this forbids the pairs $\{x,x+q/2\}$ too.

\begin{example}\label{ex:palette}
Take $q=3$ and the three sets $A_1=\{1,2\}$, $A_2=\{0,2\}$ and $A_3=\{0,1\}$ of residues modulo $3$; a set of two residues contains no progression modulo $3$, since such a progression consists of all three residues. Each residue lies in two of the sets, so its palette has $s=2$ colors: $(2,3)$ for $u=0$, $(1,3)$ for $u=1$ and $(1,2)$ for $u=2$. Starting from the valid $2$-coloring $1,1,2,2,1,1,2,2$ of $\{0,\ldots,7\}$ we obtain the valid $3$-coloring of $\{0,\ldots,23\}$ in Figure~\ref{fig:palette}. Since $W(2,3)=9$ and $W(3,3)=27$, this is close to optimal.
\end{example}

\begin{figure}[ht]
\centering
\begin{tikzpicture}[x=0.6cm,y=0.6cm]
\node[font=\scriptsize,anchor=east] at (-0.6,2.2) {$t=$};
\node[font=\scriptsize,anchor=east] at (-0.6,1.4) {short coloring of $t$};
\foreach \c/\s [count=\t from 0] in {1/k1,1/k1,2/k2,2/k2,1/k1,1/k1,2/k2,2/k2} {
 \node[\s] at (\t,1.4) {$\c$};
 \node[font=\scriptsize] at (\t,2.2) {$\t$};
}
\foreach \u/\pal/\row in {0/{(2,3)}/{2/k2,2/k2,3/k3,3/k3,2/k2,2/k2,3/k3,3/k3},1/{(1,3)}/{1/k1,1/k1,3/k3,3/k3,1/k1,1/k1,3/k3,3/k3},2/{(1,2)}/{1/k1,1/k1,2/k2,2/k2,1/k1,1/k1,2/k2,2/k2}} {
 \node[font=\scriptsize,anchor=east] at (-0.6,-\u) {$u=\u$, palette $\pal$};
 \foreach \c/\s [count=\t from 0] in \row \node[\s] at (\t,-\u) {$\c$};
}
\end{tikzpicture}
\caption{The product of Example~\ref{ex:palette}. Top: a valid $2$-coloring of $\{0,\ldots,7\}$. Below: the point $3t+u$ is in column $t$ and row $u$, and gets entry $1$ or $2$ of the palette of $u$, according to the color of $t$. Read column by column, this is a valid $3$-coloring of $\{0,\ldots,23\}$.}
\label{fig:palette}
\end{figure}

\begin{theorem}[Kim~\cite{Kim}, Theorem~2.1]\label{thm:palette}
Let $A_1,\ldots,A_r\subseteq\Z_q$ contain no progression modulo $q$, and let every residue lie in at least $s$ of them. If $\{0,\ldots,n-1\}$ has a valid $s$-coloring, then $\{0,\ldots,qn-1\}$ has a valid $r$-coloring. Hence $W(r,3)-1\geq q\bigl(W(s,3)-1\bigr)$.
\end{theorem}

\begin{proof}
For each $u\in\Z_q$ choose $s$ indices $c_1(u)<\cdots<c_s(u)$ with $u\in A_{c_i(u)}$; this is the palette of $u$. Given a valid coloring $\beta$ of $\{0,\ldots,n-1\}$ with colors $1,\ldots,s$, color $qt+u$ ($0\leq u<q$) with $c_{\beta(t)}(u)$. Suppose $x,x+d,x+2d$ all get the color $c$. If $q\nmid d$, their residues form a progression modulo $q$ inside $A_c$, which is impossible. If $d=qd'$, the three points have the same residue $u$, and $c=c_{\beta(t)}(u)$ determines $\beta(t)$; so $\beta$ is constant on the progression $t,t+d',t+2d'$, which is impossible too.
\end{proof}

Each $A_c$ has at most $r_3(\Z_q)$ elements, so the theorem needs
\begin{equation}\label{eq:counting}
 s\,q\leq\sum_{c=1}^r|A_c|\leq r\,r_3(\Z_q).
\end{equation}
Table~\ref{tab:r3} gives $r_3(\Z_q)$ for $q\leq44$.

\begin{table}[ht]
\centering
\begin{tabular}{rr@{\qquad}rr@{\qquad}rr@{\qquad}rr}
\toprule
$q$ & $r_3(\Z_q)$ & $q$ & $r_3(\Z_q)$ & $q$ & $r_3(\Z_q)$ & $q$ & $r_3(\Z_q)$\\
\midrule
2 & 1 & 13 & 4 & 24 & 6 & 35 & 9\\
3 & 2 & 14 & 4 & 25 & 7 & 36 & 8\\
4 & 2 & 15 & 4 & 26 & 7 & 37 & 10\\
5 & 2 & 16 & 4 & 27 & 8 & 38 & 8\\
6 & 2 & 17 & 5 & 28 & 8 & 39 & 10\\
7 & 3 & 18 & 5 & 29 & 8 & 40 & 9\\
8 & 3 & 19 & 6 & 30 & 8 & 41 & 9\\
9 & 4 & 20 & 5 & 31 & 8 & 42 & 9\\
10 & 4 & 21 & 6 & 32 & 8 & 43 & 9\\
11 & 4 & 22 & 6 & 33 & 8 & 44 & 9\\
12 & 4 & 23 & 6 & 34 & 8 &  & \\
\bottomrule
\end{tabular}

\caption{The largest size $r_3(\Z_q)$ of a subset of $\Z_q$ without progressions modulo $q$, computed by exhaustive search.}
\label{tab:r3}
\end{table}

\begin{proposition}\label{prop:derived}
From the bounds of Table~\ref{tab:vdw-published}, inequality~\eqref{eq:superadditive} and Theorem~\ref{thm:palette} with $q=3$ one gets $W(13,3)>2331$, $W(14,3)>2337$, $W(15,3)>2676$, $W(16,3)>2678$ and $W(17,3)>3561$.
\end{proposition}

\begin{proof}
By~\eqref{eq:superadditive}, $W(13,3)-1\geq2329+2$ and $W(14,3)-1\geq2329+8$, using $W(1,3)=3$ and $W(2,3)=9$. With $q=3$, the sets $A_c$ have at most two elements, so $s\leq\lfloor2r/3\rfloor$; this gives $3\cdot892=2676$ for $r=15,16$ ($s=10$) and $3\cdot1187=3561$ for $r=17$ ($s=11$). Finally $W(16,3)-1\geq2676+2$ by~\eqref{eq:superadditive}.
\end{proof}

For $r=13$ and $r=14$ no product reaches $2329$: by~\eqref{eq:counting} and Table~\ref{tab:r3}, the largest value of $q(W(s,3)-1)$ allowed with $q\leq44$ is $1700$ and $2289$. For $15$, $16$ and $17$ colors the counting~\eqref{eq:counting} allows more than Proposition~\ref{prop:derived}: $37\cdot75=2775$, $19\cdot170=3230$ and $27\cdot170=4590$ respectively, from $(s,q)=(4,37)$, $(5,19)$ and $(5,27)$. The question is whether the sets exist.

\subsection{Fifteen, sixteen and seventeen colors}\label{sec:vdw1516}

The progression-free subsets of $\Z_{19}$ are small: there are exactly $57$ of size $6$ and none of size $7$, and the $57$ form a single orbit of $\{0,1,3,12,14,15\}$ under the maps $x\mapsto ax+b$, $a\neq0$, which preserve progressions. A short local search over $16$-tuples of them found the cover in Table~\ref{tab:cover}: every residue modulo $19$ lies in five of its sets, except $1$, which lies in six.

\begin{table}[ht]
\centering
\begin{tabular}{llll}
\toprule
$A_{1}$&$\{4, 6, 7, 11, 12, 14\}$&$A_{9}$&$\{0, 6, 7, 10, 11, 17\}$\\
$A_{2}$&$\{4, 6, 9, 10, 13, 15\}$&$A_{10}$&$\{5, 7, 13, 14, 17, 18\}$\\
$A_{3}$&$\{3, 5, 11, 12, 15, 16\}$&$A_{11}$&$\{1, 9, 11, 14, 15, 18\}$\\
$A_{4}$&$\{0, 2, 8, 9, 12, 13\}$&$A_{12}$&$\{0, 3, 4, 7, 9, 17\}$\\
$A_{5}$&$\{1, 4, 5, 8, 10, 18\}$&$A_{13}$&$\{1, 2, 4, 13, 15, 16\}$\\
$A_{6}$&$\{0, 2, 10, 12, 15, 16\}$&$A_{14}$&$\{2, 3, 6, 8, 16, 18\}$\\
$A_{7}$&$\{1, 7, 8, 11, 12, 18\}$&$A_{15}$&$\{2, 3, 5, 14, 16, 17\}$\\
$A_{8}$&$\{0, 1, 5, 6, 8, 17\}$&$A_{16}$&$\{1, 3, 9, 10, 13, 14\}$\\
\bottomrule
\end{tabular}

\caption{Sixteen progression-free subsets of $\Z_{19}$ covering every residue at least five times.}
\label{tab:cover}
\end{table}

\begin{proposition}\label{prop:vdw16}
$W(16,3)>3230$.
\end{proposition}

\begin{proof}
Apply Theorem~\ref{thm:palette} with $q=19$, $s=5$, the sets of Table~\ref{tab:cover}, and Heule's valid $5$-coloring of $[170]$~\cite{Heule,HeuleRepo}.
\end{proof}

For $(s,q)=(4,37)$ with $15$ colors, and $(5,27)$ with $17$ colors, the same local search did not find a cover. Both are tight: by~\eqref{eq:counting}, $15\cdot10=150$ against $4\cdot37=148$, and $17\cdot8=136$ against $5\cdot27=135$. We therefore solved them exactly, as integer programs. Enlarging a set of a cover keeps it a cover, so we may assume that every $A_c$ is a maximal progression-free set; by~\eqref{eq:counting}, it then has at least $r_3(\Z_q)-\bigl(r\,r_3(\Z_q)-sq\bigr)$ elements. The program has one integer variable per such set, the number of times it is used, and asks for at most $r$ sets in total that cover every residue at least $s$ times.

For $q=37$ there are $46{,}768$ maximal progression-free sets with at least $8$ elements: $148$ of size $10$ and $46{,}620$ of size $8$ (none of size $9$ is maximal). An integer programming solver found the cover of Table~\ref{tab:cover37}. Its sets have ten elements each, only twelve of them are different, and every residue lies in four to six of them.

\begin{table}[ht]
\centering
{\small\begin{tabular}{llll}
\toprule
$A_{1}$&$\{0, 1, 4, 6, 10, 15, 17, 18, 23, 25\}$&$A_{9}$&$\{2, 3, 9, 11, 12, 14, 18, 28, 35, 36\}$\\
$A_{2}$&$\{0, 1, 6, 8, 20, 21, 24, 26, 30, 35\}$&$A_{10}$&$\{2, 4, 16, 17, 20, 22, 26, 31, 33, 34\}$\\
$A_{3}$&$\{0, 5, 7, 8, 13, 15, 27, 28, 31, 33\}$&$A_{11}$&$\{2, 4, 16, 17, 20, 22, 26, 31, 33, 34\}$\\
$A_{4}$&$\{0, 5, 7, 19, 20, 23, 25, 29, 34, 36\}$&$A_{12}$&$\{4, 5, 8, 10, 14, 19, 21, 22, 27, 29\}$\\
$A_{5}$&$\{1, 2, 5, 6, 13, 23, 27, 29, 30, 32\}$&$A_{13}$&$\{6, 7, 10, 12, 16, 21, 23, 24, 29, 31\}$\\
$A_{6}$&$\{1, 3, 15, 16, 19, 21, 25, 30, 32, 33\}$&$A_{14}$&$\{7, 8, 11, 13, 17, 22, 24, 25, 30, 32\}$\\
$A_{7}$&$\{2, 3, 9, 11, 12, 14, 18, 28, 35, 36\}$&$A_{15}$&$\{9, 10, 13, 15, 19, 24, 26, 27, 32, 34\}$\\
$A_{8}$&$\{2, 3, 9, 11, 12, 14, 18, 28, 35, 36\}$&&\\
\bottomrule
\end{tabular}
}
\caption{Fifteen progression-free subsets of $\Z_{37}$ covering every residue at least four times. Some of them coincide ($A_7=A_8=A_9$ and $A_{10}=A_{11}$), which Theorem~\ref{thm:palette} allows.}
\label{tab:cover37}
\end{table}

\begin{proposition}\label{prop:vdw15}
$W(15,3)>2775$.
\end{proposition}

\begin{proof}
Apply Theorem~\ref{thm:palette} with $q=37$, $s=4$, the sets of Table~\ref{tab:cover37}, and a valid $4$-coloring of $[75]$. Heule's certificate for four colors~\cite{HeuleRepo} colors $74$ points; one more point, of its second color, in front of them gives a valid coloring of $75$ points.
\end{proof}

For $17$ colors and $(s,q)=(5,27)$ the integer program over the $2673$ maximal progression-free subsets of $\Z_{27}$ with at least $7$ elements is infeasible. There is also a short proof.

\begin{proposition}\label{prop:z27}
No $17$ progression-free subsets of $\Z_{27}$ cover every residue five times.
\end{proposition}

\begin{proof}
Let $C_i=\{x\in\Z_{27}:x\equiv i\pmod 3\}$. The map $x\mapsto(x-i)/3$ sends the progressions inside $C_i$, whose differences are multiples of $3$, onto the progressions of $\Z_9$, so a progression-free set has at most $r_3(\Z_9)=4$ elements in each $C_i$. An exhaustive search shows that every progression-free subset of $\Z_{27}$ with $r_3(\Z_{27})=8$ elements has four elements in two of the classes $C_i$ and none in the third. Suppose that $A_1,\ldots,A_{17}$ cover every residue five times. Then $\sum_c|A_c|\geq5\cdot27=135$ and $|A_c|\leq8$, so we may number the sets so that $A_1,\ldots,A_{16}$ have $8$ elements. Let $\mu_i$ be the number of these sixteen sets that miss $C_i$. The class $C_i$ needs $5\cdot9=45$ incidences and gets at most $4(16-\mu_i)+4$, so $\mu_i\leq5$. But every one of the sixteen sets misses a class, so $\mu_0+\mu_1+\mu_2=16$, a contradiction.
\end{proof}

\subsection{Local search}\label{sec:vdw-local}

The other colorings were found by extending colorings with a local search, as for the Sidon partitions in Section~\ref{sec:repair}. A step picks a random monochromatic progression; for each of its three points it computes, for every color, how many progressions through that point would become monochromatic if the point took that color, and it makes the best move that is not tabu, with $2\%$ random moves. When no move improves, the progressions that are monochromatic at that moment get their weight increased by one, and the counts are weighted. The progressions through a point of $[N]$ are found by looping over the common difference, so one step takes $O(N)$ operations and nothing is stored except one weight per progression.

The searches started from the colorings behind Proposition~\ref{prop:derived}: for $13$ and $14$ colors, a valid $12$-coloring of $[2329]$ (ours, of the same form as Heule's) followed by $2$ or $8$ points in new colors; for $15$ and $17$ colors, the product of Theorem~\ref{thm:palette} with $q=3$ applied to Komkov's colorings for $10$ and $11$ colors; for $16$ colors, first the product with $q=3$ and then the coloring of Proposition~\ref{prop:vdw16}. Each run repeatedly adds a few points with random colors at the end and repairs them, and adds fewer points after repeated failures (Appendix~\ref{app:computations} gives the parameters and running times). The bounds in Theorem~\ref{thm:vdw} for $13$, $14$, $16$ and $17$ colors are the longest colorings these runs reached. For $15$ colors they reached only $2678$; started from the coloring of Proposition~\ref{prop:vdw15}, the same search reached $\WFifteen$ and then stopped improving. We also tried deleting the sixteenth color from prefixes of the coloring of Proposition~\ref{prop:vdw16} and repairing with $15$ colors; at lengths $2700$ to $3000$ the search remained hundreds of progressions away from a valid coloring.

\subsection{Verification}\label{sec:vdw-verify}

\begin{proof}[Proof of Theorem~\ref{thm:vdw}]
For $r=13,\ldots,17$ the repository~\cite{SRVerif} contains an $r$-coloring of $[N]$, where $N$ is the bound of Theorem~\ref{thm:vdw}. A short program in the repository checks, for each color class $C$, that no two elements $x<z$ of $C$ have $(x+z)/2\in C$, so that no class contains a progression.
\end{proof}

The same program also checks that the sets of Tables~\ref{tab:cover} and~\ref{tab:cover37} are progression-free and cover every residue five, resp.\ four, times; that the products of Propositions~\ref{prop:vdw16} and~\ref{prop:vdw15} with the colorings of $[170]$ and $[75]$ are the stored colorings of $[3230]$ and $[2775]$; the fact about $\Z_{27}$ used in Proposition~\ref{prop:z27}; and the statements above about $\Z_{19}$, and $r_3(\Z_q)$ for $q\leq20$. The values of Table~\ref{tab:r3} for $21\leq q\leq44$ come from a short exhaustive search in C++. Each coloring was also checked by a second, independent program, and the colorings of Propositions~\ref{prop:vdw15} and~\ref{prop:vdw16} by a third.

\subsection{Notes}\label{sec:vdw-notes}

Theorem~\ref{thm:palette} is due to Kim~\cite{Kim}. With $q=3$ and $A_c$ the residues other than $c\bmod3$, it is the recurrence $W(r,3)-1\geq3(W(r-\lceil r/3\rceil,3)-1)$ of Blankenship, Cummings and Taranchuk~\cite{BCT}, and with disjoint blocks of colors it is Xu's product~\cite{Xu}. The numbers $r_3(\Z_q)$ are the independence numbers of the cyclic van der Waerden hypergraphs studied by Liber~\cite{Liber}. The integer programs of Section~\ref{sec:vdw1516} were solved with HiGHS~\cite{HiGHS}, called through SciPy~\cite{SciPy}. The local search of Section~\ref{sec:vdw-local} uses a tabu list~\cite{HertzDeWerra} and the breakout method~\cite{Morris}, and the unsuccessful search of Appendix~\ref{app:computations} used CaDiCaL through PySAT~\cite{CaDiCaL,PySAT}. The structure of Heule's colorings is described in~\cite{Heule}.

\section{Open problems}\label{sec:open}

\begin{enumerate}
\item Determine $\SR(8)$ and $\SR(9)$. The bounds are now $82\leq\SR(8)\leq86$ and $103\leq\SR(9)\leq109$.
\item Can $[108]$ be partitioned into nine $12$-element Sidon sets? Such a partition would give $\SR(9)=109$. Excluding it, which by Section~\ref{sec:stop} needs an argument that uses integrality, is the first step towards $\SR(9)\leq108$.
\item The known values and our lower bounds exceed $k^2$ for $k\leq13$, equal $k^2$ for $k=14$ and $k=18$, and are slightly below $k^2$ for the other $k\leq20$ (Table~\ref{tab:bounds}). Is $\SR(k)-k^2$ bounded below?
\item How long is the longest window of Lemma~\ref{lem:window} for the Singer and Ruzsa partitions, as a function of $k$? The heuristic of Section~\ref{sec:windows} suggests order $\sqrt{k\log k}$, while Theorem~\ref{thm:general} only guarantees about $\sqrt{2k}$.
\item Heule~\cite[Conjecture~3]{Heule} conjectures $W(r,3)\geq r^3$ for all $r\geq2$. The lower bound $W(r,3)\geq2^{(c-o(1))(\log_2r)^2}$ of Section~\ref{sec:intro-notes} proves this for all sufficiently large $r$. Even ignoring the $o(1)$, that bound exceeds $r^3$ only when $\log_2r>3/c\approx21.3$, so the conjecture is really about small $r$. Our colorings exceed $r^3$ for $r=13$ but not for $14\leq r\leq17$ (Table~\ref{tab:vdw}). Are there valid colorings of $[r^3]$ for these $r$?
\end{enumerate}

\section*{Acknowledgments}

This work was supported by UNAM-PAPIIT project IN114726.

\section*{Use of AI tools}

The authors acknowledge the use of Claude (Anthropic) to assist with programming, in particular to write and optimize the search and verification programs behind the computational results, and to help with exploration, drafting, translation, and reorganization of parts of the manuscript. All mathematical content, computational output relevant to the final results, and final formulations were reviewed and verified by the authors, who approve the resulting manuscript and take full responsibility for its contents.

\begingroup
\renewcommand{\topfraction}{0.95}
\renewcommand{\bottomfraction}{0.95}
\renewcommand{\textfraction}{0.05}
\renewcommand{\floatpagefraction}{0.8}
\let\footnotesize\scriptsize
\appendix
\clearpage
\section{The integer weights}\label{app:weights}

The four columns of Table~\ref{tab:weights} are the weights of Table~\ref{tab:counts}. Zeros are included, so the weight of any set can be read directly from the table; the repository has the same vectors in \texttt{certificates.json}.
\begin{table}[htbp]
\centering\footnotesize
\setlength{\tabcolsep}{3.2pt}
\begin{tabular}{r|rrrr@{\qquad}r|rrrr@{\qquad}r|rrrr}
\toprule
$x$&$w_0$&$w_1$&$w_2$&$w'$&$x$&$w_0$&$w_1$&$w_2$&$w'$&$x$&$w_0$&$w_1$&$w_2$&$w'$\\
\midrule
1&0&1822&1941&0&38&550&0&538&0&75&464&209&366&251\\
2&2052&0&1743&0&39&435&469&0&0&76&527&404&473&0\\
3&1769&1902&0&4297&40&346&524&388&0&77&696&676&409&0\\
4&1812&1519&1985&0&41&230&357&411&483&78&523&786&945&362\\
5&2207&1806&1525&3511&42&305&338&314&0&79&742&831&964&252\\
6&1711&1963&1969&4020&43&539&558&259&142&80&420&786&586&0\\
7&1860&1652&1947&0&44&112&271&266&0&81&237&558&1057&0\\
8&0&2124&1582&0&45&198&79&138&0&82&0&324&444&0\\
9&0&0&1939&0&46&155&366&321&0&83&793&0&331&253\\
10&2089&0&0&2827&47&416&114&366&0&84&727&910&0&0\\
11&1962&1650&0&1205&48&0&82&327&0&85&771&686&1051&0\\
12&1358&1808&1803&2724&49&0&0&264&0&86&594&893&874&0\\
13&1277&1423&1776&2475&50&412&0&0&0&87&939&913&1204&0\\
14&1189&1183&1240&283&51&392&0&0&0&88&1238&1068&809&189\\
15&1274&1185&982&729&52&129&115&168&0&89&1165&1046&1046&542\\
16&1504&1088&1041&789&53&0&265&423&0&90&1093&1356&1081&0\\
17&1070&1270&951&884&54&414&0&228&0&91&874&847&1364&0\\
18&0&941&1192&0&55&0&331&117&0&92&582&811&1023&320\\
19&1237&0&823&129&56&237&217&415&0&93&966&850&836&856\\
20&1276&1071&0&1256&57&10&328&210&0&94&1347&956&1044&428\\
21&1153&1075&915&596&58&0&0&173&0&95&1275&1236&1127&410\\
22&0&1061&988&0&59&4&204&219&0&96&1163&1496&1172&1450\\
23&938&0&1082&0&60&478&0&234&0&97&1465&1350&1452&975\\
24&857&810&0&0&61&0&317&234&0&98&1196&1473&1449&950\\
25&1477&902&640&0&62&135&0&345&0&99&1434&1381&1276&2118\\
26&832&1205&824&337&63&165&93&0&0&100&1585&1652&1551&0\\
27&829&725&1102&0&64&0&145&205&0&101&1405&1554&1439&0\\
28&697&1017&732&0&65&402&0&165&0&102&0&1453&1518&0\\
29&625&956&825&0&66&394&437&0&0&103&1876&0&1516&3334\\
30&880&805&874&336&67&222&271&491&78&104&1675&1781&0&3841\\
31&756&836&460&28&68&358&421&293&0&105&0&1693&1792&0\\
32&756&956&689&0&69&133&340&443&0&106&1757&0&1553&4246\\
33&602&1006&793&0&70&0&136&322&0&107&0&1740&0&0\\
34&642&516&526&124&71&300&0&195&0&108&1834&0&1869&0\\
35&712&568&394&0&72&625&276&0&0&109&1780&1749&0&4152\\
36&597&883&532&0&73&233&586&486&287&&&&&\\
37&0&577&618&0&74&72&152&651&0&&&&&\\
\bottomrule
\end{tabular}
\caption{Integer weights for the complements of $A_0,A_1,A_2$ and $A_0\cup B_1$, respectively. Every set of the corresponding family has weight at least $10000$.}\label{tab:weights}
\end{table}

\clearpage

\section{Windows for more colors}\label{app:windows}

Table~\ref{tab:windows-large} extends Table~\ref{tab:windows} to all prime powers $q\leq64$ and all primes $p\leq67$. For each row, the repository stores a perfect difference set containing $0$ (the Singer partition of Section~\ref{sec:singer}), or the primitive root used in Section~\ref{sec:ruzsa}, together with $u$, $t$ and $m$. The program \texttt{verify\_windows.py} rebuilds each partition of $[v+m]$ and checks it, and \texttt{check\_windows\_max.py} checks, over all affine images, that $m$ is the longest window; both use only Python's standard library.
\begin{table}[htbp]
\centering
\footnotesize
\begin{tabular}{lrrrrlrr}
\toprule
Construction&$k$&$v$&$m_0$&$m$&$(u,t)$&$\SR(k)\geq$&$-k^2$\\
\midrule
Ruzsa, $p=3$&3&6&2&2&$(1,0)$&9&+0\\
Singer, $q=2$&3&7&2&3&$(1,1)$&11&+2\\
Singer, $q=3$&4&13&3&4&$(1,5)$&18&+2\\
Ruzsa, $p=5$&5&20&3&3&$(1,0)$&24&-1\\
Singer, $q=4$&5&21&3&5&$(1,17)$&27&+2\\
Singer, $q=5$&6&31&3&6&$(1,10)$&38&+2\\
Ruzsa, $p=7$&7&42&3&5&$(1,4)$&48&-1\\
Singer, $q=7$&8&57&4&8&$(2,7)$&66&+2\\
Singer, $q=8$&9&73&4&9&$(5,17)$&83&+2\\
Singer, $q=9$&10&91&4&9&$(1,86)$&101&+1\\
Ruzsa, $p=11$&11&110&4&7&$(1,3)$&118&-3\\
Singer, $q=11$&12&133&5&10&$(2,123)$&144&+0\\
Ruzsa, $p=13$&13&156&5&7&$(5,0)$&164&-5\\
Singer, $q=13$&14&183&5&12&$(17,47)$&196&+0\\
Ruzsa, $p=17$&17&272&6&8&$(1,9)$&281&-8\\
Singer, $q=16$&17&273&6&14&$(19,188)$&288&-1\\
Singer, $q=17$&18&307&6&16&$(9,293)$&324&+0\\
Ruzsa, $p=19$&19&342&6&9&$(1,9)$&352&-9\\
Singer, $q=19$&20&381&6&16&$(44,147)$&398&-2\\
Ruzsa, $p=23$&23&506&7&12&$(9,0)$&519&-10\\
Singer, $q=23$&24&553&7&19&$(4,405)$&573&-3\\
Singer, $q=25$&26&651&7&21&$(22,431)$&673&-3\\
Singer, $q=27$&28&757&7&22&$(31,418)$&780&-4\\
Ruzsa, $p=29$&29&812&7&13&$(1,4)$&826&-15\\
Singer, $q=29$&30&871&8&24&$(50,456)$&896&-4\\
Ruzsa, $p=31$&31&930&8&14&$(13,28)$&945&-16\\
Singer, $q=31$&32&993&8&23&$(89,146)$&1017&-7\\
Singer, $q=32$&33&1057&8&28&$(71,276)$&1086&-3\\
Ruzsa, $p=37$&37&1332&8&15&$(11,24)$&1348&-21\\
Singer, $q=37$&38&1407&9&27&$(29,642)$&1435&-9\\
Ruzsa, $p=41$&41&1640&9&19&$(11,39)$&1660&-21\\
Singer, $q=41$&42&1723&9&30&$(20,1423)$&1754&-10\\
Ruzsa, $p=43$&43&1806&9&16&$(1,32)$&1823&-26\\
Singer, $q=43$&44&1893&9&30&$(158,1222)$&1924&-12\\
Ruzsa, $p=47$&47&2162&10&18&$(5,11)$&2181&-28\\
Singer, $q=47$&48&2257&10&33&$(493,154)$&2291&-13\\
Singer, $q=49$&50&2451&10&33&$(20,1918)$&2485&-15\\
Ruzsa, $p=53$&53&2756&10&24&$(5,40)$&2781&-28\\
Singer, $q=53$&54&2863&10&38&$(230,2525)$&2902&-14\\
Ruzsa, $p=59$&59&3422&11&21&$(13,37)$&3444&-37\\
Singer, $q=59$&60&3541&11&39&$(110,1451)$&3581&-19\\
Ruzsa, $p=61$&61&3660&11&22&$(13,54)$&3683&-38\\
Singer, $q=61$&62&3783&11&39&$(340,2050)$&3823&-21\\
Singer, $q=64$&65&4161&11&38&$(109,1355)$&4200&-25\\
Ruzsa, $p=67$&67&4422&11&27&$(1,38)$&4450&-39\\
\bottomrule
\end{tabular}
\caption{Longest periodic windows over all affine images $x\mapsto ux+t$ of the Singer
partitions ($k=q+1$, $v=q^2+q+1$, $q\leq64$) and the Ruzsa partitions ($k=p$, $v=p(p-1)$,
$p\leq67$). Here $m_0$ is the length guaranteed by Proposition~\ref{prop:average}, $m$ is the
longest window, $(u,t)$ is one maximizer, and the last column is the lower bound $v+m+1$
minus $k^2$.}
\label{tab:windows-large}
\end{table}

\section{Computational details}\label{app:computations}

This appendix collects what is needed to repeat the computations: the programs that check the certificates, the machines and running times of the searches, and the searches that found nothing. None of it is needed to follow the arguments.

\paragraph{Checking the certificates.} All the checks are in the verification repository~\cite{SRVerif}, and each is run as \texttt{python3} followed by the command in Table~\ref{tab:commands}, from the root of the repository. They use Python's standard library; the checks of the upper bounds also need a C++17 compiler with OpenMP. All arithmetic in the checks is exact, with integers or rationals, and no solver or floating-point computation is involved.

\begin{table}[ht]
\centering\small
\begin{tabular}{@{}p{0.31\linewidth}p{0.47\linewidth}p{0.15\linewidth}@{}}
\toprule
Command&What it checks&Time\\
\midrule
\texttt{verify.py -{}-lower-only}&The partitions behind the lower bounds of Theorem~\ref{thm:main} (Section~\ref{sec:search}).&seconds\\
\texttt{verify.py -{}-threads 24}&The above, and $\SR(9)\leq109$ (Section~\ref{sec:nine}).&2 minutes$^{a}$\\
\texttt{verify\_upper.py -{}-threads 24}&$\SR(10)\leq132$, $\SR(11)\leq158$ and $\SR(12)\leq185$ (Sections~\ref{sec:upper1112} and~\ref{sec:ten}). With \texttt{-{}-g17}, also that no ruler with $17$ marks has length at most $184$; with \texttt{-{}-exhaustive}, also the exact cover search of Section~\ref{sec:computation}.&10 minutes$^{b}$; 11 and 30 with the options\\
\texttt{check\_tiling108.py}&The fractional tiling of $[108]$ (Section~\ref{sec:stop}).&under a second\\
\texttt{verify\_windows.py}&The windows of Table~\ref{tab:windows} and Appendix~\ref{app:windows}.&under a second\\
\texttt{check\_windows\_max.py}&That each of these windows is the longest one, and at least as long as Proposition~\ref{prop:average} guarantees.&40 seconds\\
\texttt{vdw/verify.py}&The colorings of Theorem~\ref{thm:vdw} and everything listed in Section~\ref{sec:vdw-verify}.&under a second\\
\texttt{verify\_examples.py}&Every small example and figure, and the other small facts stated in the text.&under a second\\
\bottomrule
\end{tabular}
\caption{The programs of the verification repository. Times: $^{a}$on a $24$-core Intel Core Ultra 9 285K; $^{b}$on a $24$-core Intel Core i9-14900KF. The values of $r_3(\Z_q)$ for $21\leq q\leq44$ in Table~\ref{tab:r3} come from the C++ program \texttt{vdw/r3\_cyclic.cpp}, which takes seconds.}
\label{tab:commands}
\end{table}

\paragraph{Machines and running times of the searches.} The exact search of Section~\ref{sec:exact} ran on one AMD Ryzen Threadripper 2990WX with $32$ cores. The exhaustive run for $[81]$ took $19.4$ CPU-hours on $16$ processes. The first partition appeared within a few minutes on eight processes for $[100]$, after about $1.5$ hours on thirteen processes for $[101]$, and after about six days on fourteen processes for $[102]$. The directory \texttt{search} of the repository has the programs, the commands of the runs and the positive controls.

The local searches of Section~\ref{sec:search} for $10$, $12$, $14$ and $19$ colors ran on the Intel Core Ultra 9 285K. Those for $9$, $11$, $13$, $15$ to $18$ and $20$ colors (for $9$ and $20$ colors, the ones that reached $[99]$ and $[396]$) ran on an Intel Core i9-14900KF and an NVIDIA RTX~4090, with one search per CPU thread or per GPU warp.

For van der Waerden numbers, the integer program for $\Z_{37}$ (Section~\ref{sec:vdw1516}) was solved in about $50$ minutes. Each run of the local search of Section~\ref{sec:vdw-local} adds $8$ points with random colors at the end and searches for $300$ seconds; after three failures in a row it adds half as many points. Each run used one thread of an AMD Ryzen Threadripper 2990WX, two runs per number of colors. For $13$, $14$, $16$ and $17$ colors the runs passed the bounds of Proposition~\ref{prop:derived} within ten minutes; for $15$ colors it took about an hour.

\paragraph{Searches that found nothing.} We record them because they may save time for others.

\emph{Sidon partitions.} Other methods were tried on the open targets of Table~\ref{tab:open}, and none found a new partition: reinserting uncolored points while every class stays Sidon, swaps that keep the class sizes fixed, and exact searches of neighborhoods. For the last, the colors of all points outside a set $U$ of free points that contains the remaining conflicts are fixed, and one decides whether $U$ can be recolored without conflicts. We decided this either with a SAT solver, a program that decides whether a system of Boolean constraints has a solution, or with a depth-first search with forward checking. Around colorings of $[288]$ with $17$ colors and of $[324]$ with $18$ colors that had one conflict, freeing the points of the conflict and up to $120$ further points usually gave an unsatisfiable formula within seconds, and with $160$ to $240$ free points the solver reached its work limit without an answer ($10^6$ conflicts in the solver's own sense, that is, dead ends of its search). Around colorings of $[224]$ with $15$ colors and of $[175]$ with $13$ colors that had one conflict, the depth-first search examined about $60\,000$ and $14\,000$ neighborhoods of $24$ to $72$ free points, and none contained a partition. This shows only that these particular colorings cannot be completed nearby.

\emph{Van der Waerden numbers with ten colors.} Many of Heule's colorings are periodic with some period $P$, and they are constant on the orbits of a group of multiplications: there is a multiplicative subgroup $H$ of the units modulo $P$ such that, measured from a center $z$, the points $z+x$ and $z+hx$ get the same color for every $h\in H$ (all modulo $P$). In his colorings for $8$ and $10$ colors the center is the largest prime factor of $P$, and the one for $W(10,3)>889$ has $P=444$ and $|H|=6$. We encoded colorings with this structure as Boolean formulas and tried every cyclic subgroup $H$ with $4\leq|H|\leq40$ for $446\leq P\leq456$, with both centers $0$ and the largest prime factor of $P$, using a SAT solver. Of $522$ formulas, $410$ were unsatisfiable and $112$ were undecided after $300$ seconds. A valid coloring of this form with $P\geq446$ would have beaten Komkov's $892$.
\clearpage
\endgroup

\bibliographystyle{plainurl}
\bibliography{references}
\end{document}